\documentclass[12pt,reqno]{amsart}
\usepackage{amsmath,amsfonts,amsthm,amsopn,amssymb,mathtools,mathrsfs,enumerate,color}
\usepackage{comment}
\usepackage{tikz}
\usetikzlibrary{arrows.meta}
\usepackage{cite,marginnote}
\usepackage{graphicx}
\usepackage[colorlinks=true,urlcolor=blue,
citecolor=red,linkcolor=blue,linktocpage,pdfpagelabels,
bookmarksnumbered,bookmarksopen]{hyperref}
\usepackage[english]{babel}
\usepackage[left=2.9cm,right=2.9cm,top=2.8cm,bottom=2.8cm]{geometry}

\numberwithin{equation}{section}

\makeindex
\newtheorem{theorem}{Theorem}[section]
\newtheorem{corollary}{Corollary}[section]
\newtheorem{remark}{Remark}[section]

\newtheorem{proposition}{Proposition}[section]
\newtheorem{lemma}{Lemma}[section]
\newtheorem{iteration lemma}{iteration Lemma}[section]

\title[Symmetry for the general-$p$ CKN inequality]{Symmetry and symmetry breaking for the general-$p$ Caffarelli--Kohn--Nirenberg inequality}

\author[S.~Fan]{Song Fan}
\author[G-D.~Li]{Gui-Dong Li}
\author[J.~J.~Zhang]{Jianjun Zhang}

\address[S.~Fan]{\newline\indent School of Mathematics and Statistics
\newline\indent Guizhou University
\newline\indent Guiyang, 550025, Guizhou, PR China}
\email{\href{mailto:doraemonsong77@gmail.com}{doraemonsong77@gmail.com}}

\address[G-D.~Li]{\newline\indent School of Mathematics and Statistics
\newline\indent Guizhou University
\newline\indent Guiyang, 550025, Guizhou, PR China}
\email{\href{mailto:bestdong123@163.com}{bestdong123@163.com}}

\address[J.~J.~Zhang]{\newline\indent College of Mathematics and Statistics
\newline\indent Chongqing Jiaotong University
\newline\indent Xuefu, Nan'an, 400074, Chongqing, PR China}
\email{\href{mailto:zhangjianjun09@tsinghua.org.cn}{zhangjianjun09@tsinghua.org.cn}}

\begin{document}

\begin{abstract}
We determine the threshold between radial symmetry and symmetry breaking
for nonnegative optimizers of the sharp Caffarelli--Kohn--Nirenberg
inequality
\[
 \left(\int_{\mathbb R^d}|x|^{-bq}|v|^q\,dx\right)^{1/q}
 \le C_{a,b}
 \left(\int_{\mathbb R^d}|x|^{-ap}|\nabla v|^p\,dx\right)^{1/p}
\]
for $d\ge2$, $1<p<d$, $a<d/p-1$, and $a<b<a+1$, with
$1/q=1/p-(1+a-b)/d$. Above and on the general-$p$ Felli--Schneider curve,
we prove that the sharp constant $C_{a,b}$ coincides with the sharp radial constant and that
every nontrivial nonnegative optimizer is radial. Combined with the linear
instability result of Caldiroli and Musina (\emph{Milan J. Math.},
\textbf{81} (2013), 421--430) below this curve, this settles their
general-$p$ symmetry conjecture, including on the threshold curve.

Unlike the case $p=2$ treated via nonlinear diffusion by Dolbeault,
Esteban, and Loss (\emph{Invent. Math.}, \textbf{206} (2016), 397--440),
our proof relies on weighted rearrangement and a $p$-independent weighted
isoperimetric inequality in dimensions $d\ge3$. We show that centered
balls uniquely minimize weighted perimeter at fixed positive volume up to
their linear stability threshold, 
confirming the two-power-weight conjecture of Alvino et al.\
(\emph{J. Math. Anal. Appl.}, \textbf{451} (2017), 280--318)
for nonnegative perimeter exponents.
Restricting this isoperimetric inequality
to radial graphs also yields a sharp inequality on $\mathbb S^{d-1}$
interpolating between the Euclidean isoperimetric inequality and the
spherical logarithmic Sobolev inequality.

\vskip0.1in


\end{abstract}

\maketitle

\section{Introduction and main results}

In variational problems invariant under a symmetry group, minimizers need
not inherit the symmetry of the functional. Identifying when symmetry
holds or breaks is both a structural question in nonlinear PDEs and a
means of computing sharp functional constants and classifying extremizers,
as in the classical Sobolev and Hardy--Littlewood--Sobolev inequalities
\cite{Aub76,Tal76,Lie83}. For a radial critical point, a negative
direction of the second variation immediately implies symmetry breaking.
Conversely, nonnegativity of the second variation is purely local: it
excludes neither a nonradial competitor with lower energy far from the
radial family nor nonradial equality cases on the threshold. 
Establishing global symmetry therefore requires global comparison or
rigidity methods, such as rearrangement \cite{Aub76,Tal76,Lie83},
tensorization or superadditivity \cite{Gro75,Car91}, or nonlinear
diffusion  \cite{DEL14,DEL16,DEL18}.

The sharp first-order Caffarelli--Kohn--Nirenberg (CKN) inequality
\cite{CKN84},
\begin{equation}
  \label{eq:ckn-norm}
 \left(\int_{\mathbb R^d}|x|^{-bq}|v|^q\,dx\right)^{1/q}
 \le C_{a,b}
 \left(\int_{\mathbb R^d}|x|^{-ap}|\nabla v|^p\,dx\right)^{1/p},
\end{equation}
provides a natural parameter-dependent model for this transition
\cite{CKN84,CW01,GM09}. We consider the admissible 
parameter range
\begin{equation}\label{eq:ckn-parameter-range}
 \begin{gathered}
 d\ge2,\qquad 1<p<d,\qquad a<a_c:=\frac dp-1,
 \qquad a\le b\le a+1,\\
 \frac1q=\frac1p-\frac{1+a-b}{d},
 \qquad p\le q\le p^*:=\frac{dp}{d-p},
 \end{gathered}
\end{equation}
where the relation between $p$ and $q$ follows from dilation invariance. Write $p'=p/(p-1)$.
On 
\[
 \mathcal D_a^{1,p}
 :=\overline{C_c^\infty(\mathbb R^d)}^{\,\|\cdot\|_{\mathcal D_a^{1,p}}},
 \qquad
 \|v\|_{\mathcal D_a^{1,p}}
 :=\left(\int_{\mathbb R^d}|x|^{-ap}|\nabla v|^p\,dx\right)^{1/p},
\]
the sharp constant $C_{a,b}$ is given by
\begin{equation}\label{eq:original-quotient}
 C_{a,b}^{-p}
 =\inf_{0\ne v\in\mathcal D_a^{1,p}}
 \frac{\displaystyle\int_{\mathbb R^d}|x|^{-ap}|\nabla v|^p\,dx}
 {\displaystyle\left(\int_{\mathbb R^d}|x|^{-bq}|v|^q\,dx\right)^{p/q}},
\end{equation}
and restricting the infimum to radial functions defines the sharp radial
constant $C_{a,b}^{\star}$. We say that \emph{symmetry is broken} if
$C_{a,b}>C_{a,b}^{\star}$; in the open strip $a<b<a+1$, we say that
\emph{symmetry holds} if $C_{a,b}=C_{a,b}^{\star}$ and every nonzero
nonnegative optimizer is radial.

For $a\le b<a+1$, introduce the effective dimension $n$ and the radial
exponent $\alpha$ by
\begin{equation}\label{eq:parameters}
 n=\frac{d}{1+a-b}\ge d,\quad
 q=\frac{np}{n-p},\quad
 \alpha=\frac{d-p-ap}{n-p}
 =\frac{(1+a-b)(a_c-a)}{a_c-a+b}>0.
\end{equation}
By Musina \cite[Theorem~1.4]{M14}, the radial infimum $(C_{a,b}^{\star})^{-p}$
is attained for $a\le b<a+1$, and its nonzero nonnegative optimizers are
precisely $A\,v_\star(\ell x)$ ($A,\ell>0$), where
\begin{equation}\label{eq:ckn-radial-profile}
 v_\star(x)=\left(1+|x|^{\alpha p'}\right)^{-(n-p)/p}.
\end{equation}
Write $d\theta$ for surface measure on $\mathbb S^{d-1}$ and set \( \omega_{d-1}:=|\mathbb S^{d-1}|\),
\( d\sigma:=\frac{d\theta}{\omega_{d-1}}\).
With $\operatorname{B}$ denoting the Euler beta function,
the sharp radial constant is given by
\cite[Theorem~2.2(2)]{HK12}
\begin{equation}\label{eq:main-constant}
 C_{a,b}^{\star}
 =\frac{\alpha^{-(n-1)/n}}{n^{1/p}}
 \left(\frac{p-1}{n-p}\right)^{(p-1)/p}
 \left[\frac{\omega_{d-1}}{p'}\,
      \operatorname{B}\!\left(\frac{n}{p'},\frac{n}{p}\right)\right]^{-1/n}.
\end{equation}

\subsection{Background and boundary lines}

The two boundary lines of \eqref{eq:ckn-parameter-range} are known. On
the Hardy line $b=a+1$ ($q=p$), one has
$C_{a,a+1}=C_{a,a+1}^{\star}=(a_c-a)^{-1}$ for all $a<a_c$, and the
infimum is unattained in $\mathcal D_a^{1,p}$
\cite[Remark~2.3(2) and Proposition~2.3]{HK12}. On the critical Sobolev
diagonal $b=a$ ($n=d$), Horiuchi and Kumlin \cite[Theorem~2.2]{HK12}
showed that $C_{a,a}=C_{a,a}^{\star}=\alpha^{-(d-1)/d}C_{0,0}$ for
$0\le a<a_c$, where $C_{0,0}$ is the sharp Euclidean Sobolev constant;
the nonzero nonnegative optimizers are $A\,v_\star(\ell x)$ for
$0<a<a_c$, and $A\,v_\star(\ell(x-x_0))$ at $a=0$ \cite{Aub76,Tal76}.
For $a<0$ on $b=a$, one has $C_{a,a}\ge C_{0,0}>C_{a,a}^{\star}$ (with
$C_{a,a}=C_{0,0}$ unattained \cite{CW01} for $p=2$).

It remains to determine the symmetry region in the open strip
$a<b<a+1$ (equivalently $p<q<p^*$ or $n>d$). For $p=2$, Chou and Chu
\cite{CC93} proved symmetry for $0\le a<a_c$, while Catrina and Wang
\cite[Theorem~1.3(ii)]{CW01} showed that for each $a<0$ optimizers are
nonradial when $b>a$ is sufficiently close to $a$; further
symmetry-breaking results for $p=2$ were obtained by Byeon and Wang
\cite{BW02}. Computing the second variation of \eqref{eq:original-quotient}
at $v_\star$ in the first spherical harmonic sector, Felli and Schneider
\cite[Corollary~1.2]{FS03} proved symmetry breaking whenever
\begin{equation}\label{eq:quadratic-fs-threshold}
 a<a_{\mathrm{FS}}(q):=a_c-2\sqrt{\frac{d-1}{q^2-4}},
 \qquad b=a-a_c+\frac{d}{q}.
\end{equation}
Following the threshold existence and bifurcation studies in
\cite{DELT09,DE12,DE14}, Dolbeault, Esteban, and Loss
\cite[Theorem~1.1]{DEL16} used nonlinear diffusion and Carr\'e du Champ
estimates on a weighted cylinder to prove that every nonnegative
optimizer is radial for $a\ge a_{\mathrm{FS}}(q)$, including on the
curve itself.

For general $1<p<d$, sharp constants, attainment, and partial symmetry
breaking were studied in \cite{H97,GM09,HK12}; in the open strip
$a<b<a+1$, the sharp constants and radial optimizers were established
for $a=0$ in \cite{H97,GY00} and for $0<a<a_c$ by Lam and Lu
\cite[Theorem~5.1]{LL17}.  Extending the linear
instability test to $1<p<d$, Caldiroli and Musina \cite[Theorem~1.1]{CM13}
proved that $v_\star$ is unstable in the first spherical harmonic mode
whenever $\alpha^2(n-1)>d-1$, giving the necessary condition
\begin{equation}\label{eq:FS}
 \alpha\le\alpha_{\mathrm{FS}}:=\sqrt{\frac{d-1}{n-1}},
\end{equation}
and conjectured that \eqref{eq:FS} is sufficient for
$C_{a,b}=C_{a,b}^{\star}$ \cite[Section~1]{CM13}. From
\eqref{eq:parameters},
\[
 \alpha^2(n-1)=\frac{(a_c-a)^2}{p^2}(q-p)\bigl((p-1)q+p\bigr),
 \qquad b=a-a_c+\frac dq.
\]
Solving $\alpha^2(n-1)=d-1$ for $b$ defines the
general-$p$ Felli--Schneider curve $b=b_{\mathrm{FS}}^{(p)}(a)$, where
\begin{equation}\label{eq:general-p-fs-curve}
 b_{\mathrm{FS}}^{(p)}(a)
 :=a-a_c+
 \frac{2d(p-1)(a_c-a)}
 {p\!\left((p-2)(a_c-a)+
 \sqrt{p^2(a_c-a)^2+4(p-1)(d-1)}\right)},
\end{equation}
so that \eqref{eq:FS} is equivalent to $b\ge b_{\mathrm{FS}}^{(p)}(a)$
and reduces to \eqref{eq:quadratic-fs-threshold} when $p=2$.
For quantitative stability of CKN inequalities, see, for example, \cite{WW22,FP24,DT25,ZZ26}

For $p\ne2$, two obstructions kept the conjecture open in dimensions
$d\ge3$ for $a<0$. First, the nonlinear diffusion method of \cite{DEL16}
relies on the $L^2$ Hilbertian structure and does not extend directly to
$p\ne2$. Ciraolo and Corso \cite{CC22} established symmetry for positive critical points in
the smaller curvature range $\alpha^2\le(d-2)/(n-2)$ under an additional
parameter restriction. Second, weighted rearrangement reduces the comparison of CKN constants
for every $1<p<d$ to a $p$-independent isoperimetric inequality with
two power weights \cite{BBMP99,ABCMP17}. In dimension $d=2$, McGillivray
proved uniqueness of the centered isoperimetric set
\cite[Theorem~1.1]{M23} and determined the threshold for the existence
of a radial CKN optimizer \cite[Theorem~14.1]{M23}. In dimensions
$d\ge3$, the weighted isoperimetric comparison and its CKN application
due to Alvino et al.\ \cite[Theorems~5.3 and~8.2]{ABCMP17} give the sufficient range
\begin{equation}\label{eq:known-geometric-symmetry}
 \alpha^2\le\alpha_0^2:=\frac{n(d-1)^2}{d(n-1)^2}
       <\frac{d-1}{n-1}=\alpha_{\mathrm{FS}}^2.
\end{equation}
Thus the intermediate regime $\alpha_0^2<\alpha^2\le\alpha_{\mathrm{FS}}^2$
and the underlying two-power-weight isoperimetric conjecture
\cite[Conjecture~5.1]{ABCMP17} remained open.

\subsection{Main results}

\begin{figure}[htbp]
\centering
\begin{tikzpicture}[x=1.45cm,y=1.1cm,font=\small,
                    >=Stealth,line cap=round,line join=round]
 \fill[black!9] (-6,-3)--(2,1)--(2,3)--(-6,-1)--cycle;
 \fill[black!30]
  (-6,-3)--(0,0)
  .. controls (-1.25,0) and (-4.5,-.68) .. (-6,-1.4)
  --cycle;
 \draw[black,line width=.85pt] (-6,-3)--(2,1);
 \draw[black,densely dashed,line width=.85pt] (-6,-1)--(2,3);
 \draw[black,line width=1.1pt]
  (-6,-1.4) .. controls (-4.5,-.68) and (-1.25,0) .. (0,0);
 \draw[black!60,densely dashed,thin] (2,-.23)--(2,3.12);
 \foreach \y in {1,3}
  \draw[black,fill=white,line width=.7pt] (2,\y) circle[radius=1.8pt];
 \draw[black,line width=.7pt] (-6.27,0)--(0,0);
 \draw[black,->,line width=.7pt] (0,0)--(2.65,0) node[right=2pt] {$a$};
 \draw[->,line width=.7pt] (0,-3.3)--(0,3.52) node[above=2pt] {$b$};
 \draw (0,2)--(-.055,2) node[left=3pt] {$1$};
 \draw (-4,0)--(-4,-.055) node[below left=3pt] {$-1$};
 \node[below=4pt] at (2,0) {$a_c$};
 \fill[black] (0,0) circle[radius=1.8pt];
 \node[below right=3pt] at (0,0) {$0$};
 \node[align=center] at (-1.18,.47)
  {Symmetry holds\\[2pt]$C_{a,b}=C_{a,b}^{\star}$};
 \node[align=center,font=\scriptsize] at (-4.77,-1.8)
  {Symmetry is broken\\[2pt]$C_{a,b}>C_{a,b}^{\star}$};
 \node at (-2.45,-.66) {$b=b_{\mathrm{FS}}^{(p)}(a)$};
 \node[rotate=21,above=3pt] at (-3.4,.3) {$b=a+1$};
 \node[rotate=21,below=4pt] at (-2.8,-1.4) {$b=a$};
 \node[anchor=north west,font=\footnotesize,text=black!70]
  at (-6.2,3.4) {Schematic: $d\ge2$, $1<p<d$};
\end{tikzpicture}
\caption{Schematic symmetry phase diagram for $a<a_c$ and $a\le b\le a+1$
(cf.\ \cite[Figure~1]{DEL16} for $p=2$). The solid curve
$b=b_{\mathrm{FS}}^{(p)}(a)$ from \eqref{eq:general-p-fs-curve} belongs to
the symmetry region; dark gray indicates the open symmetry-breaking
region, and the dashed lines mark the unattained Hardy endpoint $b=a+1$
and the excluded critical edge $a=a_c$.}
\label{fig:ckn-regions}
\end{figure}
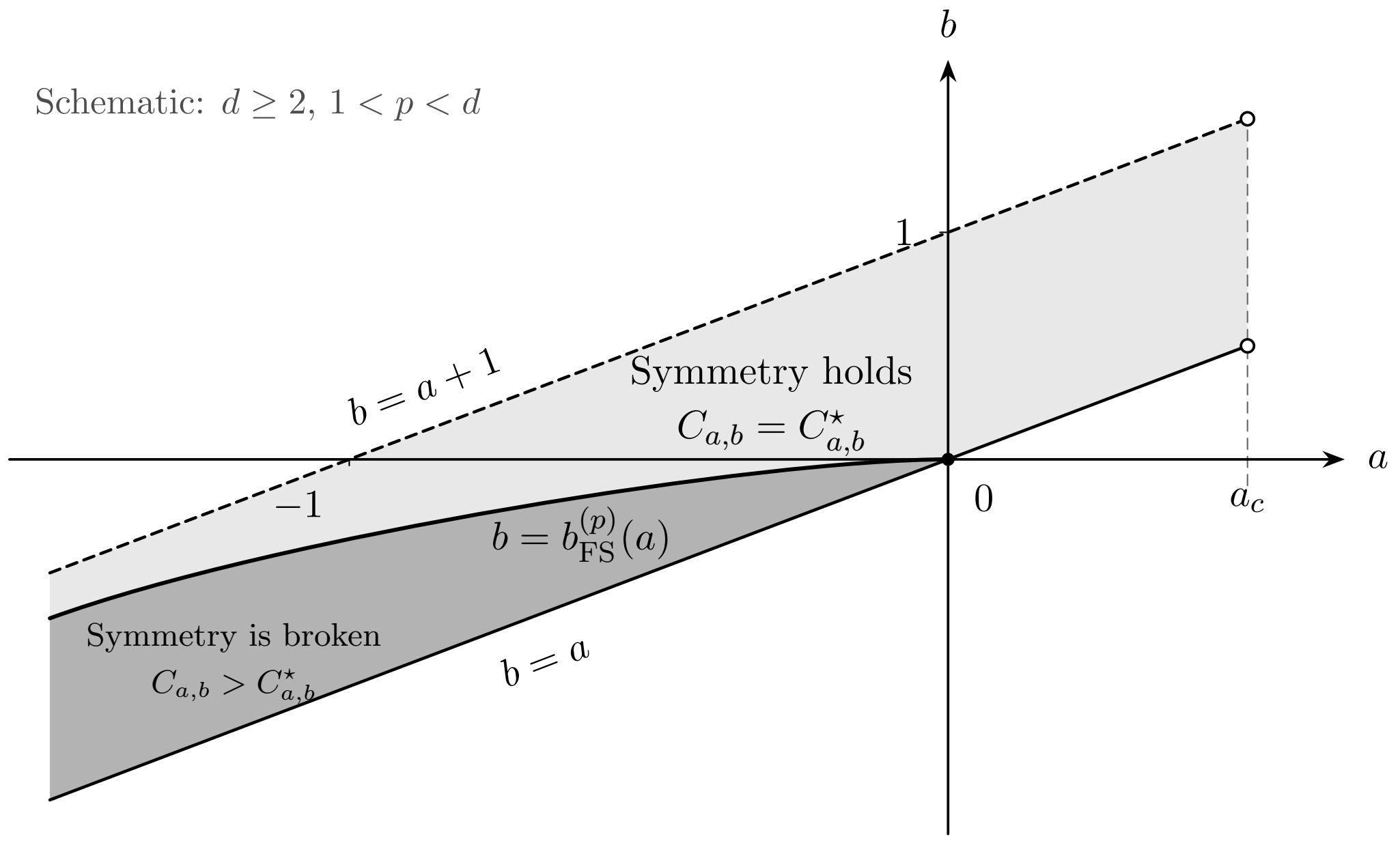

Our first main result resolves the Caldiroli--Musina conjecture
\cite{CM13} for all $d\ge2$ and $1<p<d$, including the classification of
nonnegative optimizers on $b=b_{\mathrm{FS}}^{(p)}(a)$
(Figure~\ref{fig:ckn-regions}).

\begin{theorem}[Sharp symmetry threshold for CKN]\label{thm:main-fs}
Assume \eqref{eq:ckn-parameter-range} with $a<b<a+1$. Then
\[
 \begin{array}{lll}
 b\ge b_{\mathrm{FS}}^{(p)}(a) & C_{a,b}=C_{a,b}^{\star}
       & \text{\emph{(symmetry holds)}},\\[2pt]
 b<b_{\mathrm{FS}}^{(p)}(a) & C_{a,b}>C_{a,b}^{\star}
       & \text{\emph{(symmetry is broken)}}.
 \end{array}
\]
For $b\ge b_{\mathrm{FS}}^{(p)}(a)$, $C_{a,b}$ is given by
\eqref{eq:main-constant}, and every nonzero nonnegative optimizer in
$\mathcal D_a^{1,p}$ has the form
\begin{equation}\label{eq:main-optimizers}
 v(x)=A\,v_\star(\ell x),\qquad A,\ell>0,
\end{equation}
with $v_\star$ as in \eqref{eq:ckn-radial-profile}. For
$b<b_{\mathrm{FS}}^{(p)}(a)$, every optimizer is nonradial.
\end{theorem}

In dimensions $d\ge3$, Theorem~\ref{thm:main-fs} is deduced from a
weighted isoperimetric inequality on the punctured cone
$\Omega:=(0,\infty)\times\mathbb S^{d-1}$. Identifying $(r,\theta)\in\Omega$
with $z=r\theta\in\mathbb R^d\setminus\{0\}$, define the cone metric and
weighted measure
\begin{equation}\label{eq:cone-metric-measure}
 g_\alpha=dr^2+\alpha^2r^2g_{\mathbb S^{d-1}},
 \qquad d\mu=r^{n-1}\,dr\,d\theta=r^{n-d}\,dz,
\end{equation}
where $g_{\mathbb S^{d-1}}$ is the round metric on $\mathbb S^{d-1}$. For
a set $E\subset\Omega$ of locally finite perimeter, let
$\nu_{\mathrm E}=\nu_r\theta+\nu_T$ (with $\nu_T\perp\theta$) denote its
Euclidean outer unit normal on the reduced boundary
$\partial^*E$ \cite[Chapter~15]{M12}, and let $\mathcal H^{d-1}$ denote the Euclidean
$(d-1)$-dimensional Hausdorff measure. For a Borel set $A$, the
weighted perimeter measure $P_\alpha(E;\cdot)$ is given by
\begin{equation}\label{eq:polar-perimeter}
 P_\alpha(E;A):=\int_{\partial^*E\cap\Omega\cap A}r^{n-d}
       \sqrt{\nu_r^2+\alpha^{-2}|\nu_T|^2}\,d\mathcal H^{d-1}.
\end{equation}
For $z_0\in\mathbb R^d$ and $R\ge0$, write
$B_R(z_0):=\{z\in\mathbb R^d:|z-z_0|<R\}$.
In cone integrals we restrict balls to $\Omega$ and set $\mu(\{0\})=0$.
For a positive radius function $R$ on $\mathbb S^{d-1}$, write
\begin{equation}\label{eq:radial-subgraph}
 E[R]:=\{r\theta\in\Omega:0<r<R(\theta)\}.
\end{equation}
For constant $R>0$, \eqref{eq:cone-metric-measure} and
\eqref{eq:polar-perimeter} give
\begin{equation}\label{eq:cone-centered-ball}
 \mu(B_R(0))=\frac{\omega_{d-1}}nR^n,\qquad
 P_\alpha(B_R(0);\Omega)=\omega_{d-1}R^{n-1}.
\end{equation}
For $V\ge0$, the centered ball of volume $V$ satisfies
\begin{equation}\label{eq:cone-ball-profile}
 P_\alpha(B_R(0);\Omega)=\omega_{d-1}^{1/n}(nV)^{(n-1)/n},
 \qquad R=\left(\frac{nV}{\omega_{d-1}}\right)^{1/n}.
\end{equation}

\begin{theorem}\label{thm:cone-isop}
Let $d\ge3$, $n>d$, and $0<\alpha\le\alpha_{\mathrm{FS}}$. For every
measurable set $E\subset\Omega$ of finite weighted volume and perimeter,
\begin{equation}\label{eq:cone-isop}
 P_\alpha(E;\Omega)\ge\omega_{d-1}^{1/n}\bigl(n\mu(E)\bigr)^{(n-1)/n}.
\end{equation}
If $0<\mu(E)<\infty$, equality in \eqref{eq:cone-isop} holds if and only
if $E$ coincides with a centered ball $B_R(0)$ up to a null set.
\end{theorem}

The change of variables $r=|x|^\alpha$ transforms \eqref{eq:cone-isop}
into an inequality for the weighted perimeter and volume
\[
 \operatorname{Per}_K(\widetilde E):=\int_{\partial^*\widetilde E}|x|^K\,d\mathcal H^{d-1}(x),
 \qquad
 \operatorname{Vol}_L(\widetilde E):=\int_{\widetilde E}|x|^L\,dx,
\]
confirming \cite[Conjecture~5.1]{ABCMP17} for $0\le K<L+1$ without the
convexity hypothesis assumed in \cite[Theorem~3.1]{DHKPZ19}.

\begin{corollary}\label{cor:two-power}
Let $d\ge3$, $K\ge0$, and
\begin{equation}\label{eq:two-power-range}
 0<L+1-K\le\frac{d-1}{K+d-1}.
\end{equation}
For every measurable set $\widetilde E\subset\mathbb R^d$ with finite weighted volume
and perimeter,
\begin{equation}\label{eq:two-power-isop}
 \operatorname{Per}_K(\widetilde E)\ge\omega_{d-1}^{(L+1-K)/(L+d)}
       \bigl((L+d)\operatorname{Vol}_L(\widetilde E)\bigr)^{(K+d-1)/(L+d)}.
\end{equation}
For positive volume, equality holds if and only if $\widetilde E$ coincides with a
centered ball up to a null set, except when $(K,L)=(0,0)$, where all
Euclidean balls are optimal.
\end{corollary}

\begin{remark}
When $n=d$ and $0<\alpha\le1$ (with $\alpha=1$ corresponding to
$K=L=0$),  one has
$d\mu=\alpha^{1-d}\,dV_{g_\alpha}$, and $g_\alpha$
has nonnegative Ricci curvature. Since $\alpha_{\mathrm{FS}}=1$,
\eqref{eq:cone-isop} follows by monotonicity in $\alpha$ from the
Euclidean isoperimetric inequality at $\alpha=1$ ($K=L=0$)
\cite[Theorem~14.1]{M12}. When $K=L>0$, one has $\alpha=1>\alpha_{\mathrm{FS}}$ and
centered balls are unstable \cite{RCBM08}, in contrast to log-convex
radial densities \cite{RCBM08,Cha19} (see also \cite{CRS16} for
homogeneous weights in convex cones under concavity hypotheses). For
$n>d$, a direct computation of the $n$-Bakry--\'Emery Ricci tensor
$\operatorname{Ric}_{g_\alpha,\mu}^n$ of $(\Omega,g_\alpha,d\mu)$ gives
\[
 \operatorname{Ric}_{g_\alpha,\mu}^n(\partial_r,\cdot)=0,
 \qquad
 \operatorname{Ric}_{g_\alpha,\mu}^n|_{T\mathbb S^{d-1}}
 =\bigl((d-2)-(n-2)\alpha^2\bigr)g_{\mathbb S^{d-1}}.
\]
Thus $\operatorname{Ric}_{g_\alpha,\mu}^n\ge0$ holds if and only if
$\alpha^2\le(d-2)/(n-2)$---precisely the curvature threshold reached for
critical points in \cite{CC22}---which is strictly smaller than
$\alpha_{\mathrm{FS}}^2=(d-1)/(n-1)$. In particular, curvature-dimension
comparisons do not reach $\alpha_{\mathrm{FS}}$.
\end{remark}

Evaluating Theorem~\ref{thm:cone-isop} at $\alpha=\alpha_{\mathrm{FS}}$
on smooth radial graphs $r=f(\theta)^{1/(n-1)}$ yields the following
sharp inequality on $\mathbb S^{d-1}$. We write $\nabla_\theta$ and
$\Delta_\theta$ for the gradient and the Laplace--Beltrami operator
of the round metric $g_{\mathbb S^{d-1}}$.

\begin{corollary}[Sharp inequality on $\mathbb S^{d-1}$]
\label{cor:stargraph-sphere}
Let $d\ge3$ and $n>d$. For every smooth positive function $f$ on
$\mathbb S^{d-1}$,
\begin{equation}\label{eq:stargraph-sphere}
 \int_{\mathbb S^{d-1}}
 \sqrt{f^2+\frac{|\nabla_\theta f|^2}{(d-1)(n-1)}}\,d\sigma
 \ge \left(\int_{\mathbb S^{d-1}} f^{n/(n-1)}\,d\sigma\right)^{(n-1)/n},
\end{equation}
with equality if and only if $f$ is constant.
\end{corollary}

\begin{remark}
Testing \eqref{eq:stargraph-sphere} with $f=1+\varepsilon Y_1$ for a nonzero first
spherical harmonic $Y_1$ shows that the
coefficient $1/((d-1)(n-1))$ is sharp. Moreover:
\begin{enumerate}[(i)]
 \item Setting $f=R^{n-1}$ for a smooth radius $R>0$ and letting
$n\downarrow d$ recovers the Euclidean isoperimetric inequality
$|D\mathbf1_E|(\mathbb R^d)\ge\omega_{d-1}^{1/d}(d|E|)^{(d-1)/d}$ on
$E=E[R]$ from \eqref{eq:radial-subgraph} (where $|D\mathbf1_E|(\mathbb R^d)$ and $|E|$ denote
Euclidean perimeter and Lebesgue measure), with equality for all balls
containing the origin \cite[Theorem~14.1]{M12}.

\item Setting $f=u^2$ with $\int_{\mathbb S^{d-1}}u^2\,d\sigma=1$ and
taking the right derivative at zero with respect to $1/(n-1)$ recovers the sharp spherical
logarithmic Sobolev inequality (see, e.g., \cite{MW82,BDS24})
\[
 \int_{\mathbb S^{d-1}}u^2\log u^2\,d\sigma
 \le\frac{2}{d-1}\int_{\mathbb S^{d-1}}|\nabla_\theta u|^2\,d\sigma.
\]
\end{enumerate}
\end{remark}

\subsection{Outline of the proof}

Since $P_\alpha(E;\Omega)$ is nonincreasing in $\alpha$ while $P_\alpha(B_R(0);\Omega)$ is
independent of $\alpha$, it suffices to prove Theorem~\ref{thm:cone-isop}
at $\alpha=\alpha_{\mathrm{FS}}$, where the second variation at a
centered ball vanishes on the first spherical harmonics. The proof
proceeds in three steps.
\begin{enumerate}[(1)]
\item \emph{Reduction to radial graphs (Section~\ref{sec:compactness}).}
At fixed volume $V>0$, minimizers of $P_{\alpha_{\mathrm{FS}}}(\cdot;\Omega)$
exist and have bounded density-one representatives. Spherical-cap
symmetrization about $\mathbb R e$ ($e\in\mathbb S^{d-1}$), using
\cite{CPS20}, gives a minimizer $E^{\#}$ whose density-one representative
$(E^{\#})^{(1)}$ has smooth boundary in $\Omega$.
Using volume-constrained stability, dilation and rotation Jacobi fields,
the moment identity \eqref{eq:cap-profile-moment}, and
reflection across $e^\perp$, we show that $(E^{\#})^{(1)}$ contains the
origin in its interior and that its boundary is
$\{R(\theta)\theta:\theta\in\mathbb S^{d-1}\}$ for
a smooth function $R>0$ (Proposition~\ref{prop:cap-strict-starshaped}).

\item \emph{Rigidity of stationary radial graphs (Section~\ref{sec:graph-rigidity}).}
For any smooth radial graph $r=R(\theta)>0$ that is stationary for
$P_{\alpha_{\mathrm{FS}}}(\cdot;\Omega)$ at fixed $\mu$, testing the Euler--Lagrange
equation \eqref{eq:graph-stationarity-divergence} against a
volume-preserving variation and applying Bochner's formula on
$\mathbb S^{d-1}$ reduces rigidity to the inequality
\eqref{eq:angular-bellman} (Proposition~\ref{prop:stationary-graphs}).
Adding the two identities
\eqref{eq:bellman-zero-one}--\eqref{eq:bellman-zero-two} with
gradient-dependent multipliers yields a pointwise nonnegative integrand
that vanishes only when $R$ is constant
(Proposition~\ref{pro:angular-bellman}). Hence $(E^{\#})^{(1)}=B_R(0)$ for some $R>0$,
and preservation of spherical slice measures forces
the original minimizer to coincide with $B_R(0)$ up to a null set.

\item \emph{Weighted rearrangement (Section~\ref{sec:cone-conclusion}).}
Combining Theorem~\ref{thm:cone-isop} ($d\ge3$) with the two-dimensional
isoperimetric result of \cite{M23} yields the weighted
P\'olya--Szeg\H{o} inequality \eqref{eq:weak-ps} for all $1<p<d$
(Lemma~\ref{lem:cone-rearrangement}). In the equality case, for any
nonnegative CKN optimizer $v$, almost every superlevel set $\{v>t\}$
($t>0$) must coincide with a centered ball up to a null set. Together
with the radial characterization of \cite{M14} and the instability
criterion of \cite{CM13}, this proves Theorem~\ref{thm:main-fs}.
\end{enumerate}

\section{Cone coordinates and stability of centered balls}

\subsection{The radial change of variables}

Following \cite[Section~3.1]{DEL16} ($p=2$) and \cite{CC22} (general
$p$), pass from Euclidean polar coordinates
$x=s\theta\in\mathbb R^d\setminus\{0\}$ ($s=|x|>0$,
$\theta=x/|x|\in\mathbb S^{d-1}$) to the cone coordinates
$(r,\theta)$ of \eqref{eq:cone-metric-measure} via
\begin{equation}\label{eq:cone-transform}
 \Phi(s\theta):=s^\alpha\theta,\qquad r=s^\alpha,\qquad
 w:=v\circ\Phi^{-1},
\end{equation}
with $n,\alpha$ as in \eqref{eq:parameters}. This identifies
$C_c^\infty(\mathbb R^d\setminus\{0\})$ with $C_c^\infty(\Omega)$; the
extension to $\mathcal D_a^{1,p}$ is justified in
Subsection~\ref{sec:ckn-application}.

In the radial and angular directions, the metric \eqref{eq:cone-metric-measure}
gives
\begin{equation}\label{eq:cone-gradient}
 \nabla_\alpha w
 =(\partial_r w)\partial_r+\frac{1}{\alpha^2r^2}\nabla_\theta w,\qquad
 |\nabla_\alpha w|^2
 =|\partial_r w|^2+\frac{|\nabla_\theta w|^2}{\alpha^2r^2}.
\end{equation}
The chain rule in \eqref{eq:cone-transform} and \eqref{eq:cone-gradient} gives
\begin{equation}\label{eq:cone-gradient-change}
 |\nabla v(s\theta)|
 =\sqrt{\alpha^2s^{2\alpha-2}|\partial_r w|^2+s^{-2}|\nabla_\theta w|^2}
 =\alpha s^{\alpha-1}|\nabla_\alpha w(r,\theta)|.
\end{equation}
In polar coordinates, $dr=\alpha s^{\alpha-1}\,ds$. Since
$d-p-ap=\alpha(n-p)$ and $d-bq=\alpha n$ by \eqref{eq:parameters},
\[
\begin{aligned}
 |x|^{-ap}|\nabla v|^p\,dx
  &=\alpha^{p-1}s^{d-p-ap+(p-1)\alpha}|\nabla_\alpha w|^p\,dr\,d\theta
   =\alpha^{p-1}|\nabla_\alpha w|^p\,d\mu,\\
 |x|^{-bq}|v|^q\,dx
  &=\alpha^{-1}s^{d-bq-\alpha}|w|^q\,dr\,d\theta
   =\alpha^{-1}|w|^q\,d\mu,
\end{aligned}
\]
also for noninteger $n$. Integration gives
\begin{equation}\label{eq:cone-energy-mass}
 \int_{\mathbb R^d}|x|^{-ap}|\nabla v|^p\,dx
 =\alpha^{p-1}\int_{\Omega}|\nabla_\alpha w|^p\,d\mu,
 \qquad
 \int_{\mathbb R^d}|x|^{-bq}|v|^q\,dx
 =\alpha^{-1}\int_{\Omega}|w|^q\,d\mu.
\end{equation}

For $w\in W^{1,p}_{\mathrm{loc}}(\Omega)$ satisfying
$\int_{\Omega}|\nabla_\alpha w|^p\,d\mu<\infty$ and
$0<\int_{\Omega}|w|^q\,d\mu<\infty$,
\begin{equation}\label{eq:cone-quotient-factor}
 \frac{\int_{\mathbb R^d}|x|^{-ap}|\nabla v|^p\,dx}
      {\bigl(\int_{\mathbb R^d}|x|^{-bq}|v|^q\,dx\bigr)^{p/q}}
 =\alpha^{p-1+p/q}\mathcal Q_{\alpha,p}[w],\qquad
 \mathcal Q_{\alpha,p}[w]:=
 \frac{\int_{\Omega}|\nabla_\alpha w|^p\,d\mu}
      {\bigl(\int_{\Omega}|w|^q\,d\mu\bigr)^{p/q}}.
\end{equation}
For fixed $n,p$, \eqref{eq:cone-gradient} gives
$\mathcal Q_{\alpha,p}[w]\ge\mathcal Q_{\alpha_{\mathrm{FS}},p}[w]$
whenever $0<\alpha\le\alpha_{\mathrm{FS}}$, with equality for radial $w$.

\subsection{Pullback, weighted coarea, and the second variation of a ball}

For $0\le s<R<\infty$, write
\begin{equation}\label{eq:cone-annuli}
 A_{s,R}:=\{y\in\mathbb R^d:s<|y|<R\}.
\end{equation}
Use $|x|$ in the original coordinates and $|z|$ on the cone. The
exponents in the pullback by \eqref{eq:cone-transform} are
\begin{equation}\label{eq:cone-two-power-densities}
 K:=\alpha(n-1)-d+1,\qquad L:=\alpha n-d.
\end{equation}
At $x=s\theta$, the differential $d\Phi_x$ in \eqref{eq:cone-transform}
scales $\theta$ by $\alpha s^{\alpha-1}$ and $\theta^\perp$ by
$s^{\alpha-1}$. Hence $\det(d\Phi_x)=\alpha s^{d(\alpha-1)}$, and the
cofactor matrix $\operatorname{cof}(d\Phi_x)=(\det d\Phi_x)(d\Phi_x)^{-T}$
acts on $\tilde\nu_{\mathrm E}=\tilde\nu_r\theta+\tilde\nu_T$
($\tilde\nu_T\perp\theta$) by
\[
 \operatorname{cof}(d\Phi_x)(\tilde\nu_r\theta+\tilde\nu_T)
 =s^{(d-1)(\alpha-1)}(\tilde\nu_r\theta+\alpha\tilde\nu_T).
\]
Let $\widetilde E=\Phi^{-1}(E)$ and take its Euclidean outer unit normal.
Since $\alpha(n-d)+(d-1)(\alpha-1)=K$, the area formula for sets
of locally finite perimeter \cite{M12} and \eqref{eq:polar-perimeter} give
\begin{equation}\label{eq:cone-pullback-volume-perimeter}
 \mu(E)=\alpha\int_{\widetilde E}|x|^L\,dx,\qquad
 P_\alpha(E;A)=\int_{\partial^*\widetilde E\cap\Phi^{-1}(A\cap\Omega)}
               |x|^K\,d\mathcal H^{d-1}(x).
\end{equation}
For $u\in W^{1,1}_{\mathrm{loc}}(\Omega)$ and $\tilde u=u\circ\Phi$,
\eqref{eq:cone-gradient-change} gives
\[
 |\nabla_\alpha u(z)|\,d\mu(z)
 =(\alpha^{-1}s^{1-\alpha}|\nabla_x\tilde u(x)|)\,(\alpha s^{\alpha n-d}\,dx)
 =|x|^K|\nabla_x\tilde u(x)|\,dx.
\]
For each $t\in\mathbb R$,
\begin{equation}\label{eq:cone-superlevels}
 \Phi^{-1}\bigl(\{z\in\Omega:u(z)>t\}\bigr)
 =\{x\in\mathbb R^d\setminus\{0\}:\tilde u(x)>t\}.
\end{equation}
On compact annuli we may apply the Euclidean coarea formula
\cite[Theorem~13.1]{M12} with weight $|x|^K$. By \eqref{eq:cone-superlevels}
and exhaustion, every nonnegative Borel function $\Psi$ on $\Omega$ satisfies
\begin{equation}\label{eq:cone-coarea}
 \int_\Omega \Psi\,|\nabla_\alpha u|\,d\mu
   =\int_{\mathbb R}\int_{\partial^*\{u>t\}\cap\Omega}
       \Psi(z)\,P_\alpha(\{u>t\};dz)\,dt.
\end{equation}

To verify that $\alpha_{\mathrm{FS}}$ is the linear stability threshold
for centered balls (cf.\ \cite[proof of Theorem~4.1]{ABCMP17}), let
$u\in C^\infty(\mathbb S^{d-1})$ satisfy
$\int_{\mathbb S^{d-1}}u\,d\sigma=0$ and set
\begin{equation}\label{eq:ball-variation-radius}
 R_\varepsilon(\theta):=
 \frac{e^{\varepsilon u(\theta)}}
      {\left(\int_{\mathbb S^{d-1}}e^{n\varepsilon u}\,d\sigma\right)^{1/n}}.
\end{equation}
By \eqref{eq:radial-subgraph}, $\mu(E[R_\varepsilon])=\mu(B_1(0))$.
For any smooth $R>0$, the Euclidean gradient of $r-R(\theta)$ on the
graph is $\theta-\nabla_\theta\log R$. Hence
\begin{equation}\label{eq:graph-normal-and-area}
 \nu_{\mathrm E}
 =\frac{\theta-\nabla_\theta\log R}{\sqrt{1+|\nabla_\theta\log R|^2}},
 \qquad
 d\mathcal H^{d-1}
 =R^{d-1}\sqrt{1+|\nabla_\theta\log R|^2}\,d\theta.
\end{equation}
By \eqref{eq:graph-normal-and-area} and \eqref{eq:polar-perimeter},
\begin{equation}\label{eq:radial-graph-perimeter}
 P_\alpha(E[R];\Omega)
 =\int_{\mathbb S^{d-1}} R^{n-1}
       \sqrt{1+\alpha^{-2}|\nabla_\theta\log R|^2}\,d\theta.
\end{equation}
For $R_\varepsilon$ in \eqref{eq:ball-variation-radius},
$\nabla_\theta\log R_\varepsilon=\varepsilon\nabla_\theta u$.
Since $u$ has mean zero, \eqref{eq:radial-graph-perimeter} gives
\[
\begin{aligned}
 \frac{P_\alpha(E[R_\varepsilon];\Omega)}{P_\alpha(B_1(0);\Omega)}
 &=\frac{\displaystyle\int_{\mathbb S^{d-1}}e^{(n-1)\varepsilon u}
       \sqrt{1+\alpha^{-2}\varepsilon^2|\nabla_\theta u|^2}\,d\sigma}
      {\displaystyle\left(\int_{\mathbb S^{d-1}}e^{n\varepsilon u}\,d\sigma\right)^{(n-1)/n}}\\
 &=\frac{1+\frac{\varepsilon^2}{2}\int_{\mathbb S^{d-1}}
         \bigl(\alpha^{-2}|\nabla_\theta u|^2+(n-1)^2u^2\bigr)\,d\sigma+O(\varepsilon^3)}
        {1+\frac{n(n-1)\varepsilon^2}{2}\int_{\mathbb S^{d-1}}u^2\,d\sigma+O(\varepsilon^3)}\\
 &=1+\frac{\varepsilon^2}{2}\int_{\mathbb S^{d-1}}
   \bigl(\alpha^{-2}|\nabla_\theta u|^2-(n-1)u^2\bigr)\,d\sigma
   +O(\varepsilon^3).
\end{aligned}
\]
The first variation vanishes. The first positive eigenvalue of
$-\Delta_\theta$ is $d-1$, so the second variation is nonnegative
on mean-zero functions exactly when $\alpha^2(n-1)\le d-1$.
At $\alpha=\alpha_{\mathrm{FS}}$ its kernel is the first spherical
harmonic space. For $\alpha>\alpha_{\mathrm{FS}}$ it is negative on
each nonzero element of that space, for example $u(\theta)=e\cdot\theta$
with $e\in\mathbb S^{d-1}$.

\section{Reduction of isoperimetric minimizers to radial graphs}
\label{sec:compactness}

Throughout this section we assume $n>d\ge3$.

\subsection{Spherical-cap symmetrization}
Write $X=X_r\theta+X_T$, $X_T\perp\theta$, for
$X\in C_c^1(\Omega;\mathbb R^d)$ at $z=r\theta$. Then
$|X|_{g_\alpha}^2=X_r^2+\alpha^2|X_T|^2$, and Cauchy--Schwarz gives
\[
 X\cdot\nu_{\mathrm E}
 =X_r\nu_r+(\alpha X_T)\cdot(\alpha^{-1}\nu_T)
 \le\sqrt{X_r^2+\alpha^2|X_T|^2}\sqrt{\nu_r^2+\alpha^{-2}|\nu_T|^2},
\]
with equality for positively proportional $(X_r,\alpha X_T)$ and
$(\nu_r,\alpha^{-1}\nu_T)$. Gauss--Green on compact annuli and
exhaustion of $\Omega$ give the dual formula (cf.\ \cite[Chapters~12 and~15]{M12})
\begin{equation}\label{eq:cap-dual-perimeter}
 P_\alpha(E;\Omega)=\sup\left\{\int_{E\cap\Omega}\operatorname{div}_z(r^{n-d}X)\,dz
 : X\in C_c^1(\Omega;\mathbb R^d),\ X_r^2+\alpha^2|X_T|^2\le1\right\},
\end{equation}
and hence lower semicontinuity under $L^1_{\mathrm{loc}}(\Omega)$
convergence. The weighted divergence and Laplacian are
\begin{equation}\label{eq:cone-weighted-operators}
 \operatorname{div}_\mu X:=r^{d-n}\operatorname{div}_z(r^{n-d}X),
 \qquad \Delta_\mu u:=\operatorname{div}_\mu(\nabla_\alpha u).
\end{equation}

Fix $e\in\mathbb S^{d-1}$ and set
\begin{equation}\label{eq:cap-polar-angle}
 \vartheta(\theta):=\arccos(\theta\cdot e),\qquad
 G_e:=\{O\in SO(d):Oe=e\}.
\end{equation}
For a measurable set $E\subset\Omega$, write
\begin{equation}\label{eq:cap-spherical-slices}
 E_r:=\{\theta\in\mathbb S^{d-1}:r\theta\in E\},\qquad
 |E_r|:=\int_{E_r}d\theta .
\end{equation}
For a.e.\ $r>0$, define $\phi_E(r)\in[0,\pi]$ by
\begin{equation}\label{eq:cap-opening-angle}
 |E_r|=|\mathbb S^{d-2}|\int_0^{\phi_E(r)}\sin^{d-2}\!u\,du .
\end{equation}
For a.e.\ $r>0$, the \emph{spherical-cap symmetral} $E^{\#}$ has slice
\[
 (E^{\#})_r=
 \begin{cases}
  \varnothing, & \phi_E(r)=0,\\
  \{\theta\in\mathbb S^{d-1}:\theta\cdot e>\cos\phi_E(r)\},
      & 0<\phi_E(r)<\pi,\\
  \mathbb S^{d-1}, & \phi_E(r)=\pi.
 \end{cases}
\]
Thus $(E^{\#})_r$ is centered at $e$ and has the same spherical
measure as $E_r$ \cite{CPS20}. We call $E$ \emph{cap-symmetric about $e$}
if $\mu(E\triangle E^{\#})=0$.

\begin{lemma}\label{lem:cone-cap-symmetrization}
Let $\alpha>0$. For every measurable set $E\subset\Omega$ of finite
weighted volume and perimeter,
\[
 \mu(E^{\#})=\mu(E),\qquad P_\alpha(E^{\#};\Omega)\le P_\alpha(E;\Omega).
\]
\end{lemma}

\begin{proof}
The map \eqref{eq:cone-transform} preserves $\theta$, so the pullback
$\widetilde E$ from \eqref{eq:cone-pullback-volume-perimeter} satisfies
\[
 (\widetilde E)_s=E_{s^\alpha},\qquad
 ((\widetilde E)^{\#})_s=(E^{\#})_{s^\alpha}
 \qquad\text{for a.e. }s>0.
\]
Thus $(\widetilde E)^{\#}=\Phi^{-1}(E^{\#})$ a.e. By preservation of
slice measures,
\[
 \mu(E^{\#})
 =\int_0^\infty r^{n-1}|(E^{\#})_r|\,dr
 =\int_0^\infty r^{n-1}|E_r|\,dr
 =\mu(E).
\]

On $A_{1/j,j}$ from \eqref{eq:cone-annuli}, for integers $j\ge2$, the bound
$|x|^K\ge j^{-|K|}$ and \eqref{eq:cone-pullback-volume-perimeter} give
\[
 |D\mathbf1_{\widetilde E}|(A_{1/j,j})
 \le j^{|K|}P_\alpha(E;\Omega)<\infty.
\]
The truncation $\widetilde E\cap A_{1/j,j}$ has finite Euclidean
volume and perimeter. Its cap symmetral is $(\widetilde E)^{\#}\cap
A_{1/j,j}$ a.e., so \cite[Theorem~1.1, (1.4)]{CPS20} gives locally
finite perimeter for $(\widetilde E)^{\#}$ in $\Omega$. For Borel
$\Lambda\subset(1/j,j)$, truncation leaves both perimeter measures
unchanged inside the open annulus. Hence \cite[(1.4)]{CPS20} gives
\begin{equation}\label{eq:cap-annular-perimeter}
\begin{aligned}
 |D\mathbf1_{(\widetilde E)^{\#}}|(\{|x|\in\Lambda\})
 &=|D\mathbf1_{(\widetilde E)^{\#}\cap A_{1/j,j}}|(\{|x|\in\Lambda\})\\
 &\le |D\mathbf1_{\widetilde E\cap A_{1/j,j}}|(\{|x|\in\Lambda\})
 =|D\mathbf1_{\widetilde E}|(\{|x|\in\Lambda\}).
\end{aligned}
\end{equation}
For Borel $\Lambda\subset(0,\infty)$, \eqref{eq:cap-annular-perimeter}
on $\Lambda\cap(1/j,j)$ and continuity from below give
\begin{equation}\label{eq:cap-radial-perimeter}
 |D\mathbf1_{(\widetilde E)^{\#}}|(\{|x|\in \Lambda\})
 \le |D\mathbf1_{\widetilde E}|(\{|x|\in \Lambda\}).
\end{equation}

Integrating \eqref{eq:cap-radial-perimeter} for
$\Lambda=\{s>0:s^K>t\}$, we obtain by Tonelli's theorem and
\eqref{eq:cone-pullback-volume-perimeter}
\[
\begin{aligned}
 P_\alpha(E^{\#};\Omega)
 &=\int_\Omega |x|^K\,d|D\mathbf1_{(\widetilde E)^{\#}}|(x)
 =\int_0^\infty |D\mathbf1_{(\widetilde E)^{\#}}|(\{|x|^K>t\})\,dt\\
 &\le\int_0^\infty |D\mathbf1_{\widetilde E}|(\{|x|^K>t\})\,dt
=\int_\Omega |x|^K\,d|D\mathbf1_{\widetilde E}|(x)
  =P_\alpha(E;\Omega).
\end{aligned}
\]
This completes the proof.
\end{proof}

\subsection{Existence, boundedness, and axial regularity}

Since $n>d$, the density $r^{n-d}$ grows at infinity. The following
exterior estimate quantifies the resulting tightness of minimizing
sequences (cf.\ \cite{MP13,BDKS18}).

\begin{lemma}\label{lem:cone-exterior-volume}
Let $0<\alpha\le1$. There exists $C_d>0$, depending only on $d$, such
that for every measurable $E\subset\Omega$ of finite weighted volume and
perimeter and every $R>0$,
\begin{equation}\label{eq:exterior-perimeter-tightness}
 \mu(E\setminus B_R(0))^{(d-1)/d}
 \le C_d R^{-(n-d)/d}\,P_\alpha(E;\{r>R\}).
\end{equation}
\end{lemma}

\begin{proof}
Since $n>d$ and $\alpha\le1$, $E\cap\{r>R\}$ has finite  volume
 and perimeter in $\{r>R\}$. For $t>0$,
set $s_t:=\max\{R,t^{1/(n-d)}\}$, so that
$\{r>s_t\}=\{r>R\}\cap\{r^{n-d}>t\}$ is the exterior of the closed
Euclidean ball $\overline{B_{s_t}(0)}$. By the exterior isoperimetric
inequality \cite{CGR07} (in its finite-perimeter form
\cite[Theorem~1.1]{FusMor23}),
\[
\begin{aligned}
 |E\cap\{r>s_t\}|
 &=|E\cap\{r>s_t\}|^{1/d}|E\cap\{r>s_t\}|^{(d-1)/d}\\
 &\le C_d\,|E\cap\{r>R\}|^{1/d}\,|D\mathbf1_E|(\{r>s_t\}).
\end{aligned}
\]
Since $d\mu=r^{n-d}\,dz$ and $\alpha\le1$ (whence
$r^{n-d}|D\mathbf1_E|\le P_\alpha(E;\cdot)$ as measures), Tonelli's theorem gives
\begin{align*}
 \mu(E\setminus B_R(0))
 &=\int_0^\infty |E\cap\{r>s_t\}|\,dt\\
 &\le C_d\,|E\cap\{r>R\}|^{1/d}
       \int_0^\infty |D\mathbf1_E|(\{r>s_t\})\,dt\\
 &=C_d\,|E\cap\{r>R\}|^{1/d}
       \int_{\{r>R\}} r^{n-d}\,d|D\mathbf1_E|\\
 &\le C_d\,R^{-(n-d)/d}\,
       \mu(E\setminus B_R(0))^{1/d}\,P_\alpha(E;\{r>R\}).
\end{align*}
Cancellation gives \eqref{eq:exterior-perimeter-tightness} whenever
$\mu(E\setminus B_R(0))>0$.
\end{proof}

For $V>0$, define
\begin{equation}\label{eq:cone-isoperimetric-profile}
 I_\alpha(V):=\inf\{P_\alpha(E;\Omega):E\subset\Omega,\ \mu(E)=V\}.
\end{equation}
The density-one representative of a measurable $E\subset\Omega$ is
\begin{equation}\label{eq:cap-density-one}
 E^{(1)}:=\left\{z\in\mathbb R^d:
   \lim_{\varepsilon\downarrow0}
   \frac{\mu(E\cap B_\varepsilon(z))}{\mu(B_\varepsilon(z))}=1\right\}.
\end{equation}

\begin{lemma}\label{lem:geometric-compactness}
Let $0<\alpha\le1$ and $V>0$. Then $I_\alpha(V)$ is attained, and for
every minimizer $E$ the representative $E^{(1)}$ is bounded. If
$E^{(1)}$ is invariant under $G_e$ from \eqref{eq:cap-polar-angle}, then
$\partial E^{(1)}\cap\Omega$ is a smooth embedded hypersurface.
\end{lemma}

  \begin{proof}
\emph{Existence.}
Testing with a centered ball gives $I_\alpha(V)<\infty$.
 Let  $\mu(E_j)=V$ and $P_\alpha(E_j;\Omega)\to I_\alpha(V)$.
By \eqref{eq:cone-annuli} and
Lemma~\ref{lem:cone-exterior-volume},
\begin{equation}
  \label{eq:tightness-estimate}
\begin{aligned}
 \sup_j\mu\bigl(E_j\setminus A_{\varepsilon,R}\bigr)
 &\le \mu(B_\varepsilon(0))+\sup_j\mu(E_j\setminus B_R(0))\\
 &\le \frac{\omega_{d-1}}{n}\varepsilon^n
    +\Bigl(C_d\sup_j P_\alpha(E_j;\Omega)\Bigr)^{d/(d-1)}
     R^{-(n-d)/(d-1)}\\
 &\to 0  \qquad(\varepsilon\downarrow0,\ R\uparrow\infty).
\end{aligned}
\end{equation}
Since $\alpha\le1$, \eqref{eq:polar-perimeter} gives
$\sup_j |D\mathbf1_{E_j}|(A_{\varepsilon,R})\le\varepsilon^{d-n}\sup_j P_\alpha(E_j;\Omega)<\infty$.
BV compactness on annuli \cite[Chapter~12]{M12} and
\eqref{eq:tightness-estimate} therefore give, along a subsequence,
$\mathbf1_{E_j}\to\mathbf1_E$ in $L^1(\Omega,d\mu)$. Hence $\mu(E)=V$,
and lower semicontinuity from \eqref{eq:cap-dual-perimeter} yields
\[
 I_\alpha(V)\le P_\alpha(E;\Omega)
 \le\liminf_{j\to\infty}P_\alpha(E_j;\Omega)=I_\alpha(V).
\]

\emph{Boundedness.}
Let $E$ be any minimizer of volume $V$, and set
\[
 V_{\mathrm{out}}(R):=\mu(E\setminus B_R(0))=\int_R^\infty s^{n-1}|E_s|\,ds.
\]
For a.e.\ $R>0$,
$P_\alpha(E\cap B_R(0);\partial B_R(0))=R^{n-1}|E_R|=-V_{\mathrm{out}}'(R)$.
If $V_{\mathrm{out}}(R)<V$, dilating $E\cap B_R(0)$ by
$(V/(V-V_{\mathrm{out}}(R)))^{1/n}$ restores volume $V$. Minimality gives,
for a.e.\ such $R$,
\begin{equation}\label{eq:compactness-cut}
 \left(1-\frac{V_{\mathrm{out}}(R)}V\right)^{(n-1)/n}P_\alpha(E;\Omega)
 \le P_\alpha(E\cap B_R(0);\Omega)
 =P_\alpha(E;\Omega)-P_\alpha(E;\{r>R\})-V_{\mathrm{out}}'(R).
\end{equation}
For a.e.\ $R\ge1$ with $V_{\mathrm{out}}(R)<V$,
 \eqref{eq:exterior-perimeter-tightness}, \eqref{eq:compactness-cut},
and $1-(1-t)^{(n-1)/n}\le t$ for $0\le t\le1$ gives
\begin{equation}\label{eq:cap-truncation}
 C_d^{-1}V_{\mathrm{out}}(R)^{(d-1)/d}
 \le P_\alpha(E;\{r>R\})
 \le -V_{\mathrm{out}}'(R)+\frac{P_\alpha(E;\Omega)}V\,V_{\mathrm{out}}(R).
\end{equation}
Since $V_{\mathrm{out}}(R)\downarrow0$, choose $R_0\ge1$ with $V_{\mathrm{out}}(R_0)<V$
and $(P_\alpha(E;\Omega)/V)V_{\mathrm{out}}(R_0)^{1/d}\le 1/(2C_d)$.
Then \eqref{eq:cap-truncation} implies
\begin{equation}\label{eq:cap-tail-differential-inequality}
 \bigl(V_{\mathrm{out}}^{1/d}\bigr)'(R)\le-\frac{1}{2dC_d}
 \qquad\text{for a.e. }R\ge R_0\text{ with }V_{\mathrm{out}}(R)>0.
\end{equation}
Integrating \eqref{eq:cap-tail-differential-inequality} gives
$V_{\mathrm{out}}(R_1)=0$ for $R_1=R_0+2dC_d V_{\mathrm{out}}(R_0)^{1/d}$.
Hence $E$ has zero weighted volume outside $\overline{B_{R_1}(0)}$, and
$E^{(1)}\subset\overline{B_{R_1}(0)}$.

\emph{Axial regularity.}
Suppose $E^{(1)}$ is $G_e$-invariant. The density
$d\mu/dV_{g_\alpha}=\alpha^{1-d}r^{n-d}$ is smooth and positive on
$\Omega$. By \cite[Proposition~3.5, Corollary~3.8 and
Section~3.10]{M03}, $\partial E^{(1)}\cap\Omega$ is smooth and embedded
outside a relatively closed singular set $\mathcal S$ with
$\mathcal H^{d-2}(\mathcal S)=0$. Since $G_e\mathcal S=\mathcal S$,
\[
 z\in\mathcal S\setminus\mathbb Re
 \quad\Longrightarrow\quad
 0<\mathcal H^{d-2}(G_e z)
 \le\mathcal H^{d-2}(\mathcal S)=0.
\]
Hence $\mathcal S\subset\mathbb R e$.

Fix $z\in\partial E^{(1)}\cap\mathbb R e\cap\Omega$. In orthonormal
normal coordinates at $z$, let $E_{\mathrm{tan}}$ be the Euclidean
density-one representative of a tangent cone to $E^{(1)}$. It minimizes
Euclidean perimeter and is invariant under rotations fixing $e$
(cf.\ \cite[proof of Theorem~3.2, Steps three and four]{FigMag13}).
The orbit argument gives smoothness off the axis. Each component of
the off-axis link is a latitude, whose polar angle
$\vartheta_0\in(0,\pi)$ satisfies, by stationarity,
\[
 H_{\partial E_{\mathrm{tan}}}(t,\vartheta_0)
 =\pm\frac{d-2}{t}\cot\vartheta_0=0
 \quad\Longrightarrow\quad \vartheta_0=\frac{\pi}{2}.
\]
Since $\mathcal H^{d-1}(\mathbb R e)=0$, the perimeter of
$E_{\mathrm{tan}}$ is concentrated on $e^\perp$. Thus its indicator is
constant a.e.\ on each open half-space. Positive density of the tangent
boundary at the origin gives $|D\mathbf1_{E_{\mathrm{tan}}}|(B_1(0))>0$
\cite[Proposition~3.2 and Section~3.10]{M03}. The two constants therefore
differ, so $E_{\mathrm{tan}}$ is a half-space up to a null set.
By \cite[Proposition~3.5]{M03}, $z$ is regular. Thus $\mathcal S=\varnothing$.
\end{proof}

\subsection{Strict star-shapedness of a cap-symmetric minimizer}

For the rest of Section~\ref{sec:compactness}, fix
$\alpha=\alpha_{\mathrm{FS}}$. By the Lebesgue density theorem,
\eqref{eq:cap-density-one}, and polar integration,
\begin{equation}\label{eq:cap-representative-slices}
 \mu(E\triangle E^{(1)})=0,\qquad
 \bigl|E_r\triangle(E^{(1)})_r\bigr|=0\quad\text{for a.e. }r>0.
\end{equation}
For cap-symmetric $E$, \eqref{eq:cap-polar-angle}--\eqref{eq:cap-opening-angle}
give $\mu(OE\triangle E)=0$ for every $O\in G_e$. Such rotations preserve
$\mu$ and map balls to balls, so \eqref{eq:cap-density-one} gives
$O(E^{(1)})=E^{(1)}$.

Let $E$ be a cap-symmetric minimizer. By Lemma~\ref{lem:geometric-compactness},
$\Sigma:=\partial E^{(1)}\cap\Omega$ is smooth and embedded. Its
intersection with a fixed closed half-plane bounded by $\mathbb R e$
is the \emph{meridian}. With $\nu$ the outer unit normal in
$g_{\alpha_{\mathrm{FS}}}$, the dilation flow and its normal speed at $t=0$ are
\begin{equation}\label{eq:cap-dilation-speed}
 D_t(r,\theta):=(e^t r,\theta),\qquad
 \psi:=\langle r\partial_r,\nu\rangle_{g_{\alpha_{\mathrm{FS}}}}.
\end{equation}
For a smooth defining function $F$ on an open $U\subset\Omega$ with
$E^{(1)}\cap U=\{F<0\}$, $\Sigma\cap U=\{F=0\}$, and $dF\ne0$ on $\Sigma\cap U$,
\eqref{eq:cap-dilation-speed} gives
\begin{equation}\label{eq:cap-defining-function-speed}
 \psi
 =\frac{\langle r\partial_r,\nabla_{\alpha_{\mathrm{FS}}}F
         \rangle_{g_{\alpha_{\mathrm{FS}}}}}
        {|\nabla_{\alpha_{\mathrm{FS}}}F|}
 =\frac{r\,\partial_r F}
        {|\nabla_{\alpha_{\mathrm{FS}}}F|}
 \qquad\text{on }\Sigma\cap U.
\end{equation}

\begin{proposition}\label{prop:cap-strict-starshaped}
Let $V>0$, and let $E$ attain $I_{\alpha_{\mathrm{FS}}}(V)$ and be
cap-symmetric about $e$. Then $0\in\operatorname{int}(E^{(1)})$ and $\psi>0$
on $\Sigma=\partial E^{(1)}\cap\Omega$. In particular, $\partial_r F>0$ on
$\Sigma\cap U$ for every local defining function $F$ as in
\eqref{eq:cap-defining-function-speed}. Each ray from the origin meets
$\partial E^{(1)}$ once, and
\[
 \partial E^{(1)}=\{R(\theta)\theta:\theta\in\mathbb S^{d-1}\}
 \qquad\text{for a unique }R\in C^\infty(\mathbb S^{d-1}),\quad R>0.
\]
\end{proposition}

In Proposition~\ref{prop:cap-strict-starshaped}, volume-constrained
stability gives connectedness of $\{\psi>0\}$, in the spirit of
Barbosa and do Carmo \cite{BdC84} and the connectedness criterion for
strictly log-concave densities \cite[Corollary~3.9]{RCBM08}. Rotational
Jacobi fields give the graph $\vartheta=\phi(r)$ for the off-axis
meridian in the nonball case. For parallel ambient fields see
\cite[Lemma~3.1]{R21}, and for meridians of cap-symmetric minimizers
see \cite[Section~2]{Cha19}.

\medskip
\noindent\emph{Variational equations and stability on $\Sigma$.}
Write $dA_\mu:=\alpha_{\mathrm{FS}}^{1-d}r^{n-d}\,dA_{g_{\alpha_{\mathrm{FS}}}}$
for the weighted area element on $\Sigma$. Under \eqref{eq:cap-dilation-speed},
\begin{equation}\label{eq:cap-dilation-laws}
 \mu(D_tE)=e^{nt}V,\qquad
 P_{\alpha_{\mathrm{FS}}}(D_tE;\Omega)
 =e^{(n-1)t}P_{\alpha_{\mathrm{FS}}}(E;\Omega).
\end{equation}
By \eqref{eq:cone-isoperimetric-profile} and \eqref{eq:cap-dilation-laws},
$I_{\alpha_{\mathrm{FS}}}(s)
 =I_{\alpha_{\mathrm{FS}}}(V)(s/V)^{(n-1)/n}$ for $s>0$.
Since $E$ minimizes perimeter at volume $V$,
\begin{equation}\label{eq:cap-profile-derivatives}
 \lambda_E:=I_{\alpha_{\mathrm{FS}}}'(V)
 =\frac{(n-1)P_{\alpha_{\mathrm{FS}}}(E;\Omega)}{nV}>0,\qquad
 I_{\alpha_{\mathrm{FS}}}''(V)=-\frac{\lambda_E}{nV}.
\end{equation}

Let $s\mapsto E(s)$ be a smooth variation compactly supported in $\Omega$,
with $E(0)=E^{(1)}$ and outer normal speed $f\in C_c^\infty(\Sigma)$. Then
\begin{equation}\label{eq:cap-profile-contact}
 P_{\alpha_{\mathrm{FS}}}(E(s);\Omega)
       -I_{\alpha_{\mathrm{FS}}}(\mu(E(s)))
 \ge P_{\alpha_{\mathrm{FS}}}(E;\Omega)-I_{\alpha_{\mathrm{FS}}}(V)=0.
\end{equation}
By \eqref{eq:cap-profile-derivatives} and \eqref{eq:cap-profile-contact},
the first variation gives
\[
 0=\int_\Sigma(H_\mu-\lambda_E)f\,dA_\mu
 \qquad\text{for all }f\in C_c^\infty(\Sigma).
\]
Thus
\begin{equation}\label{eq:cap-constant-mean-curvature}
 H_\mu:=\operatorname{div}_\Sigma\nu
 +\frac{n-d}{r}\langle\nu,\partial_r\rangle_{g_{\alpha_{\mathrm{FS}}}}
 =\lambda_E\qquad\text{on }\Sigma.
\end{equation}

With $\dot H_\mu=-L_\Sigma f$, the weighted Jacobi operator and index
form of $P_{\alpha_{\mathrm{FS}}}(\cdot;\Omega)-\lambda_E\mu(\cdot)$ at
$E^{(1)}$ take the form (cf.\ \cite[Lemma~2.3, (2.6)--(2.7) and (2.9)--(2.11)]{R21})
\begin{equation}\label{eq:cap-jacobi-operator}
 \Delta_{\Sigma,\mu} f:=r^{d-n}\operatorname{div}_\Sigma(r^{n-d}\nabla_\Sigma f),
 \qquad L_\Sigma f=\Delta_{\Sigma,\mu}f+U_\Sigma f,
\end{equation}
\begin{equation} \label{eq:cap-index-form}
   Q_\Sigma(f)
 =-\int_\Sigma fL_\Sigma f\,dA_\mu
  =\int_\Sigma\bigl(|\nabla_\Sigma f|^2-U_\Sigma f^2\bigr)\,dA_\mu,
\end{equation}
where $U_\Sigma\in C^\infty(\Sigma)$. Since
$\left.\frac{d}{ds}\right|_{s=0}\mu(E(s))=\int_\Sigma f\,dA_\mu$,
the second variation of \eqref{eq:cap-profile-contact} gives, by
\eqref{eq:cap-profile-derivatives} and \eqref{eq:cap-index-form},
\begin{equation}\label{eq:cap-profile-second-variation}
 0\le Q_\Sigma(f)+\frac{\lambda_E}{nV}
   \left(\int_\Sigma f\,dA_\mu\right)^2.
\end{equation}
By smooth approximation on a common compact support,
\eqref{eq:cap-profile-second-variation} holds for every compactly
supported $f\in H^1(\Sigma)$. In particular,
\begin{equation}\label{eq:cap-mean-zero-stability}
 \int_\Sigma f\,dA_\mu=0\quad\Longrightarrow\quad Q_\Sigma(f)\ge0.
\end{equation}

Since $D_t^*g_{\alpha_{\mathrm{FS}}}=e^{2t}g_{\alpha_{\mathrm{FS}}}$
and $r\circ D_t=e^tr$, \eqref{eq:cap-constant-mean-curvature} gives
\begin{equation}\label{eq:cap-dilated-curvature}
 H_\mu[D_t\Sigma](D_tz)=e^{-t}H_\mu[\Sigma](z)=e^{-t}\lambda_E
 \qquad(z\in\Sigma).
\end{equation}
By \eqref{eq:cap-constant-mean-curvature}, $\nabla_\Sigma H_\mu=0$,
so only the normal speed in \eqref{eq:cap-dilation-speed} contributes
to the derivative of \eqref{eq:cap-dilated-curvature} at $t=0$. Thus
\begin{equation}\label{eq:cap-dilation-jacobi}
 -L_\Sigma\psi
 =\left.\frac{d}{dt}\right|_{t=0}
   \bigl(H_\mu[D_t\Sigma]\circ D_t\bigr)
 =-\lambda_E,
 \qquad L_\Sigma\psi=\lambda_E\quad\text{on }\Sigma.
\end{equation}

\begin{lemma}\label{lem:cap-positive-connected}
Under the hypotheses of Proposition~\ref{prop:cap-strict-starshaped},
$\Omega_+:=\{z\in\Sigma:\psi(z)>0\}$ is nonempty and connected, and
$\Sigma$ is connected.
\end{lemma}

\begin{proof}
Since $\Omega$ is connected and $0<\mu(E)<\mu(\Omega)$, we have
$\Sigma\ne\varnothing$. Each connected component $\Sigma_0$ of $\Sigma$
is bounded and closed in $\Omega$, so $r$ attains its maximum $R_0>0$
at some $p\in\Sigma_0$. There $\nabla_\Sigma r=0$ and
$\nu=\epsilon_p\partial_r$ for some $\epsilon_p\in\{-1,1\}$. By
$\nabla^2_{g_{\alpha_{\mathrm{FS}}}}r
=r^{-1}(g_{\alpha_{\mathrm{FS}}}-dr\otimes dr)$ and
\eqref{eq:cap-constant-mean-curvature},
\[
 0\ge\Delta_\Sigma r(p)
 =\frac{d-1}{R_0}-\epsilon_p\operatorname{div}_\Sigma\nu(p)
 =\frac{d-1}{R_0}
   -\epsilon_p\left(\lambda_E-\frac{n-d}{R_0}\epsilon_p\right)
 =\frac{n-1}{R_0}-\epsilon_p\lambda_E.
\]
As $\lambda_E>0$, we have $\epsilon_p=1$ and $\psi(p)=R_0>0$. Thus
\begin{equation}\label{eq:cap-positive-meets-components}
 \Omega_+\cap\Sigma_0\ne\varnothing
 \qquad\text{for every connected component $\Sigma_0$ of $\Sigma$.}
\end{equation}

Let $D$ be a component of $\Omega_+$ and $\psi_D:=\mathbf1_D\psi$.
Since $D$ is open in $\Sigma$ and closed in $\Omega_+$, $\psi$ vanishes
on $\partial_\Sigma D$. Truncating at height $t>0$ gives
\begin{equation}\label{eq:cap-component-gradient}
 \mathbf1_D(\psi-t)_+\to\psi_D
 \quad\text{in }H^1_{\mathrm{loc}}(\Sigma)\quad(t\downarrow0),
 \qquad
 \nabla_\Sigma\psi_D=\mathbf1_D\nabla_\Sigma\psi
 \quad\text{a.e. on }\Sigma.
\end{equation}

Choose a smooth radial cutoff $0\le\chi_\varepsilon\le1$ satisfying
\begin{equation}\label{eq:cone-radial-cutoff}
 \chi_\varepsilon(r)=0\quad(r\le\varepsilon),\qquad
 \chi_\varepsilon(r)=1\quad(r\ge2\varepsilon),\qquad
 |\chi'_\varepsilon(r)|\le C/\varepsilon .
\end{equation}
Since $\Sigma$ is bounded and $\chi_\varepsilon$ vanishes near the
vertex, $\chi_\varepsilon^2\psi_D$ has compact support in $\Sigma$.
Testing \eqref{eq:cap-dilation-jacobi} with $\chi_\varepsilon^2\psi_D$
gives, by \eqref{eq:cap-index-form} and \eqref{eq:cap-component-gradient},
\begin{equation}\label{eq:cap-component-index}
 Q_\Sigma(\chi_\varepsilon\psi_D)
 =-\lambda_E\int_D\chi_\varepsilon^2\psi\,dA_\mu
   +\int_D\psi^2|\nabla_\Sigma\chi_\varepsilon|^2\,dA_\mu.
\end{equation}
Since $|\psi|\le r$ and $\Sigma$ has finite weighted area,
\[
 0\le\int_D\psi^2|\nabla_\Sigma\chi_\varepsilon|^2\,dA_\mu
 \le 4C^2\int_{\Sigma\cap A_{\varepsilon,2\varepsilon}}dA_\mu
 \to 0.
\]
Thus \eqref{eq:cap-component-index} yields
\begin{equation}\label{eq:cap-component-negative}
 \lim_{\varepsilon\downarrow0}Q_\Sigma(\chi_\varepsilon\psi_D)
 =-\lambda_E\int_D\psi\,dA_\mu<0.
\end{equation}

If $D_1,D_2$ are distinct components, \eqref{eq:cap-component-negative}
gives $Q_\Sigma(f_i)<0$ for $f_i:=\chi_\varepsilon\psi_{D_i}$,
$i=1,2$, and $\varepsilon$ small. Put $m_i:=\int_\Sigma f_i\,dA_\mu>0$
and $f:=m_2f_1-m_1f_2$. By disjointness of $D_1,D_2$,
\eqref{eq:cap-component-gradient}, and \eqref{eq:cap-index-form},
\[
 \int_\Sigma f\,dA_\mu=0,\qquad
 Q_\Sigma(f)=m_2^2Q_\Sigma(f_1)+m_1^2Q_\Sigma(f_2)<0,
\]
contrary to \eqref{eq:cap-mean-zero-stability}. Thus $\Omega_+$ has at
most one component and, by \eqref{eq:cap-positive-meets-components},
meets every component of the nonempty $\Sigma$. Both sets are connected.
\end{proof}

Set $r_-:=\inf_{z\in\Sigma}|z|$ and $r_+:=\sup_{z\in\Sigma}|z|$.
Since $\Sigma$ is nonempty and bounded, $0\le r_-\le r_+<\infty$.

\begin{lemma}\label{lem:cap-meridian}
Under the hypotheses of Proposition~\ref{prop:cap-strict-starshaped},
either $E^{(1)}=B_{(n-1)/\lambda_E}(0)$, or $r_-<r_+$ and there exists a
unique $\phi\in C^\infty((r_-,r_+);(0,\pi))$ such that, for every
$r\in(r_-,r_+)$,
\begin{equation}\label{eq:cap-meridian-slices}
 (E^{(1)})_r=\{\theta\in\mathbb S^{d-1}:\vartheta(\theta)<\phi(r)\},
 \quad
 \Sigma\cap\partial B_r(0)
 =\{r\theta:\theta\in\mathbb S^{d-1},\ \vartheta(\theta)=\phi(r)\}.
\end{equation}
For each positive endpoint $r_0\in\{r_-,r_+\}$,
\begin{equation}\label{eq:cap-meridian-endpoints}
 \ell_0:=\lim_{\substack{r\to r_0\\r\in(r_-,r_+)}}\phi(r)
 \in\{0,\pi\},
 \qquad
 \Sigma\cap\partial B_{r_0}(0)=\{r_0\cos(\ell_0)e\}.
\end{equation}
\end{lemma}

\begin{proof}
Set $\Sigma^\circ:=\Sigma\setminus\mathbb R e$. At each axial point $p$,
rotational invariance gives $\nu(p)\parallel e$ and $T_p\Sigma=e^\perp$.
The local graph over $e^\perp$ meets the axis only at $p$, so the
axial points form a closed discrete subset of $\Sigma$. Removing these
points from the connected manifold $\Sigma$
(Lemma~\ref{lem:cap-positive-connected}) of dimension $d-1\ge2$ gives
\begin{equation}\label{eq:cap-off-axis-connected}
 \Sigma^\circ\ \text{is connected},\qquad
 \overline{\Sigma^\circ}^{\,\Omega}=\Sigma.
\end{equation}
Write $z=r(\cos\vartheta\,e+\sin\vartheta\,\upsilon)\in\Sigma^\circ$,
with $\upsilon\in\mathbb S^{d-1}\cap e^\perp$, and put
$e_\vartheta:=(\alpha_{\mathrm{FS}}r)^{-1}\partial_\vartheta$.
Axial symmetry places $\nu$ in the meridian plane. Orient the unit
tangent $\tau$ so that $(\tau,\nu)$ has the orientation of
$(\partial_r,e_\vartheta)$. Since $\partial_\tau r=dr(\tau)$ and
$dr(\nu)=\psi/r$ by \eqref{eq:cap-dilation-speed},
\begin{equation}\label{eq:cap-meridian-frame}
 \tau=(\partial_\tau r)\partial_r-\frac{\psi}{r}e_\vartheta,\qquad
 \nu=\frac{\psi}{r}\partial_r+(\partial_\tau r)e_\vartheta,\qquad
 (\partial_\tau r)^2+\frac{\psi^2}{r^2}=1 .
\end{equation}
For $0\le f\in C_c^\infty(\Omega\setminus\mathbb R e)$, Gauss--Green
and integration over the caps of opening $\phi_E(r)$ in
\eqref{eq:cap-opening-angle} give
\begin{equation}\label{eq:cap-angular-flux}
\begin{aligned}
 \int_{\Sigma^\circ}\alpha_{\mathrm{FS}}r\,f\,(\partial_\tau r)\,dA_\mu
 &=\int_E\operatorname{div}_\mu(f\,\partial_\vartheta)\,d\mu\\
 &=\int_0^\infty r^{n-1}\sin^{d-2}\!\bigl(\phi_E(r)\bigr)
   \int_{\mathbb S^{d-2}}f(r,\phi_E(r),\upsilon)\,d\upsilon\,dr
 \ge0.
\end{aligned}
\end{equation}
By continuity, \eqref{eq:cap-angular-flux} gives $\partial_\tau r\ge0$
on $\Sigma^\circ$.

For each unit $\upsilon_0\perp e$, the rotations generated by
$K_{\upsilon_0}(z):=(z\cdot e)\upsilon_0-(z\cdot\upsilon_0)e$
preserve $g_{\alpha_{\mathrm{FS}}}$ and $\mu$. The linearization
of \eqref{eq:cap-constant-mean-curvature} along these rotations gives
\cite[Lemma~3.3(i) and Remark~3.4]{Rosales23}
\begin{equation}\label{eq:cap-rotation-jacobi}
 \langle K_{\upsilon_0},\nu\rangle_{g_{\alpha_{\mathrm{FS}}}}
 =\alpha_{\mathrm{FS}}r(\upsilon\cdot\upsilon_0)\,\partial_\tau r,\qquad
 L_\Sigma\bigl(\langle K_{\upsilon_0},\nu\rangle_{g_{\alpha_{\mathrm{FS}}}}\bigr)
 =0\qquad\text{on }\Sigma^\circ.
\end{equation}
At a zero $z_0$ of $\partial_\tau r$, choose $\upsilon_0=\upsilon(z_0)$.
Since $\upsilon\cdot\upsilon_0>0$ near $z_0$, the Jacobi field in
\eqref{eq:cap-rotation-jacobi} is nonnegative there and vanishes at $z_0$.
Harnack's inequality forces it, and hence $\partial_\tau r$, to vanish
near $z_0$. The zero set is thus open and closed in $\Sigma^\circ$.
By \eqref{eq:cap-off-axis-connected},
\[
 \partial_\tau r\equiv0\quad\text{on }\Sigma^\circ
 \qquad\text{or}\qquad
 \partial_\tau r>0\quad\text{on }\Sigma^\circ.
\]
In the first case,
$|\nabla_\Sigma r|^2=1-\langle\nu,\partial_r\rangle_{g_{\alpha_{\mathrm{FS}}}}^2
=(\partial_\tau r)^2=0$ on $\Sigma^\circ$. By
\eqref{eq:cap-off-axis-connected}, $r\equiv R>0$ on $\Sigma$. As an
embedded hypersurface, $\Sigma$ is open and closed in $\partial B_R(0)$,
hence equals that sphere. Boundedness and positive volume give
$E^{(1)}=B_R(0)$, with $R=(n-1)/\lambda_E$ by
\eqref{eq:cap-constant-mean-curvature}.

In the second case, every defining function $F$ as in
\eqref{eq:cap-defining-function-speed} satisfies
\begin{equation}\label{eq:cap-angular-transversality}
 \partial_\vartheta F
 =\alpha_{\mathrm{FS}}r\,
   |\nabla_{\alpha_{\mathrm{FS}}}F|
   \partial_\tau r>0
 \qquad\text{on }\Sigma^\circ\cap U.
\end{equation}
On a fixed meridian the implicit function theorem gives a smooth local
graph $\vartheta=\phi_J(r)$, $r\in J$. Its rotations form a latitude.
By \eqref{eq:cap-angular-transversality}, the smaller-angle side lies
in $E^{(1)}$ and the larger-angle side in its exterior. By cap symmetry
and \eqref{eq:cap-representative-slices}, $\phi_J$ is the cap opening
for a.e.\ $r\in J$. Continuity on overlapping intervals gives
\begin{equation}\label{eq:cap-meridian-overlap}
 \phi_J=\phi_E\quad\text{a.e. on }J,\qquad
 \phi_{J_1}=\phi_{J_2}\quad\text{on }J_1\cap J_2,
\end{equation}
The local graphs make $r(\Sigma^\circ)$ open, and
\eqref{eq:cap-off-axis-connected} gives $r(\Sigma^\circ)=(r_-,r_+)$.
By \eqref{eq:cap-meridian-overlap} the graphs define a unique smooth
$\phi$ on this interval. In $A_{r_-,r_+}$, $E$ agrees a.e.\ with
$\{\vartheta<\phi(r)\}$. Taking density-one representatives in
\eqref{eq:cap-density-one} gives \eqref{eq:cap-meridian-slices} at
every $r\in(r_-,r_+)$.

At a positive endpoint $r_0$, compactness on closed annuli,
\eqref{eq:cap-off-axis-connected}, and \eqref{eq:cap-meridian-slices} give
\begin{equation}\label{eq:cap-endpoint-cluster}
 \bigcap_{\delta>0}
 \overline{\phi\bigl((r_-,r_+)\cap(r_0-\delta,r_0+\delta)\bigr)}
 =
 \{\vartheta(z/|z|):z\in\Sigma\cap\partial B_{r_0}(0)\}
 \subset\{0,\pi\}.
\end{equation}
Since $r_0\notin r(\Sigma^\circ)$, all cluster angles in
\eqref{eq:cap-endpoint-cluster} belong to $\{0,\pi\}$. Their set is a
nested intersection of nonempty compact intervals, hence a singleton.
This proves \eqref{eq:cap-meridian-endpoints}.
\end{proof}

\begin{lemma}\label{lem:cone-meridian-system}
In the nonball case of Lemma~\ref{lem:cap-meridian}, orient $\tau$
toward increasing $r$. Rotational invariance makes $\psi$ a function of
$r\in(r_-,r_+)$ on $\Sigma$. With primes denoting $r$-derivatives,
\begin{equation}\label{eq:cap-meridian-components}
 \partial_\tau r
 =\frac{1}{\sqrt{1+\alpha_{\mathrm{FS}}^2r^2\phi'^2}}>0,\qquad
 \psi=-\alpha_{\mathrm{FS}}r^2\phi'\,\partial_\tau r,
\end{equation}
and $(\phi,\psi)$ satisfies
\begin{align}
 r\phi'&=-\frac{\psi}{\alpha_{\mathrm{FS}}\sqrt{r^2-\psi^2}},
 \label{eq:cap-angle-ode}\\
 r\psi'&=\lambda_E r^2-(n-2)\psi
 -\frac{d-2}{\alpha_{\mathrm{FS}}}\cot\phi\,\sqrt{r^2-\psi^2}.
 \label{eq:cap-slope-ode}
\end{align}
\end{lemma}

\begin{proof}
Along $\vartheta=\phi(r)$ the unit tangent toward increasing $r$ is
\[
 \tau=\frac{\partial_r+\alpha_{\mathrm{FS}}r\phi'e_\vartheta}
             {\sqrt{1+\alpha_{\mathrm{FS}}^2r^2\phi'^2}}.
\]
Its radial component gives the first identity in
\eqref{eq:cap-meridian-components}. By \eqref{eq:cap-dilation-speed}
and \eqref{eq:cap-meridian-frame},
\begin{equation}\label{eq:cap-meridian-speed}
 \frac{\psi}{r}=-\langle\tau,e_\vartheta\rangle_{g_{\alpha_{\mathrm{FS}}}}
 =-\alpha_{\mathrm{FS}}r\phi'\,\partial_\tau r,\qquad
 \partial_\tau r=\sqrt{1-\psi^2/r^2}.
\end{equation}
Thus $|\psi|<r$ and \eqref{eq:cap-angle-ode} follows.

In coordinates $(r,\vartheta,\upsilon)$, $\upsilon\in\mathbb S^{d-2}$,
the volume and graph area elements are
\begin{equation}\label{eq:cap-meridian-measures}
\begin{aligned}
 d\mu
 &=r^{n-1}\sin^{d-2}\!\vartheta\,dr\,d\vartheta\,d\upsilon,\\
 dA_\mu
 &=\alpha_{\mathrm{FS}}^{1-d}r^{n-d}
   (\alpha_{\mathrm{FS}}r\sin\phi)^{d-2}
   \frac{dr}{\partial_\tau r}\,d\upsilon
 =\frac{r^{n-2}\sin^{d-2}\!\phi}{\alpha_{\mathrm{FS}}\,\partial_\tau r}
   \,dr\,d\upsilon .
\end{aligned}
\end{equation}
For $r_-<s<R<r_+$, equations \eqref{eq:cone-annuli},
\eqref{eq:cap-meridian-measures}, and \eqref{eq:cap-meridian-components} give
\[
 \frac{P_{\alpha_{\mathrm{FS}}}(E;A_{s,R})-\lambda_E\mu(E\cap A_{s,R})}
      {|\mathbb S^{d-2}|}
 =\int_s^R\mathcal L(r,\phi,\phi')\,dr,
\]
where
\begin{equation}\label{eq:cap-meridian-lagrangian}
 \mathcal L(r,\phi,\phi')
 =\frac{r^{n-2}\sin^{d-2}\!\phi}{\alpha_{\mathrm{FS}}}
   \sqrt{1+\alpha_{\mathrm{FS}}^2r^2\phi'^2}
  -\lambda_E r^{n-1}\int_0^\phi \sin^{d-2}\!t\,dt.
\end{equation}
By \eqref{eq:cap-meridian-lagrangian} and \eqref{eq:cap-meridian-components},
$\partial_{\phi'}\mathcal L=-r^{n-2}\sin^{d-2}\!\phi\,\psi$. For
variations $\phi+\varepsilon u$, $u\in C_c^\infty((s,R))$,
the Euler--Lagrange equation reads
\begin{equation}\label{eq:cap-meridian-euler}
 -\bigl(r^{n-2}\sin^{d-2}\!\phi\,\psi\bigr)'
 =\frac{(d-2)r^{n-2}\sin^{d-3}\!\phi\cos\phi}
        {\alpha_{\mathrm{FS}}\,\partial_\tau r}
  -\lambda_E r^{n-1}\sin^{d-2}\!\phi .
\end{equation}
Dividing \eqref{eq:cap-meridian-euler} by $r^{n-2}\sin^{d-2}\phi>0$ yields
\[
 \psi'+\frac{n-2}{r}\psi
 =\lambda_E r-(d-2)\cot\phi
   \left(\frac{1}{\alpha_{\mathrm{FS}}\,\partial_\tau r}+\phi'\psi\right).
\]
By \eqref{eq:cap-meridian-speed}, the factor in parentheses is
\[
 \frac{1-\psi^2/r^2}{\alpha_{\mathrm{FS}}\,\partial_\tau r}
 =\frac{\sqrt{r^2-\psi^2}}{\alpha_{\mathrm{FS}}r}.
\]
Multiplication by $r$ gives \eqref{eq:cap-slope-ode}.
\end{proof}

\begin{lemma}\label{lem:cap-angular-moment}
In the nonball case of Lemma~\ref{lem:cap-meridian},
\begin{equation}\label{eq:cap-profile-moment}
 \int_{r_-}^{r_+}r^{n-3}\sin^{d-2}\!\phi\,(\partial_\tau r)\cos\phi\,dr=0.
\end{equation}
In particular, $\phi(r_0)>\pi/2$ for some $r_0\in(r_-,r_+)$.
\end{lemma}

\begin{proof}
For $Z_e=r\cos\vartheta$, \eqref{eq:cone-weighted-operators} and
$\alpha_{\mathrm{FS}}^2(n-1)=d-1$ give
\begin{equation}\label{eq:cap-axial-derivatives}
\begin{aligned}
 \Delta_\mu Z_e
 &=\frac{\cos\vartheta}{r}
   \left(n-1-\frac{d-1}{\alpha_{\mathrm{FS}}^2}\right)=0,\\
 \nabla^2_{g_{\alpha_{\mathrm{FS}}}}Z_e
 &=(1-\alpha_{\mathrm{FS}}^{-2})\frac{\cos\vartheta}{r}
   (g_{\alpha_{\mathrm{FS}}}-dr\otimes dr).
\end{aligned}
\end{equation}
By \eqref{eq:cap-meridian-frame},
$(g_{\alpha_{\mathrm{FS}}}-dr\otimes dr)(\nu,\nu)=(\partial_\tau r)^2$.
The weighted tangential trace of \eqref{eq:cap-axial-derivatives}, using
\eqref{eq:cap-jacobi-operator} and $H_\mu=\lambda_E$ from
\eqref{eq:cap-constant-mean-curvature}, is
\begin{equation}\label{eq:cap-moment-trace}
\begin{aligned}
 \Delta_{\Sigma,\mu}Z_e
 &=\Delta_\mu Z_e-\nabla^2_{g_{\alpha_{\mathrm{FS}}}}Z_e(\nu,\nu)
   -\lambda_E\,\partial_\nu Z_e\\
 &=(\alpha_{\mathrm{FS}}^{-2}-1)
   \frac{\cos\phi}{r}(\partial_\tau r)^2-\lambda_E\,\partial_\nu Z_e.
\end{aligned}
\end{equation}
Along the meridian, $\partial_\tau=(\partial_\tau r)\,d/dr$ on
functions of $r$. By \eqref{eq:cap-meridian-measures} and
\eqref{eq:cap-moment-trace},
\begin{equation}\label{eq:cap-moment-product}
\begin{aligned}
 \bigl(r^{n-2}\sin^{d-2}\!\phi\,\partial_\tau Z_e\bigr)'
 &=\frac{r^{n-2}\sin^{d-2}\!\phi}{\partial_\tau r}\Delta_{\Sigma,\mu}Z_e\\
 &=(\alpha_{\mathrm{FS}}^{-2}-1)r^{n-3}\sin^{d-2}\!\phi\,(\partial_\tau r)\cos\phi
   -\lambda_E\frac{r^{n-2}\sin^{d-2}\!\phi}{\partial_\tau r}\partial_\nu Z_e.
\end{aligned}
\end{equation}
By \eqref{eq:cap-meridian-frame},
$\partial_\nu Z_e=(\psi/r)\cos\phi-\alpha_{\mathrm{FS}}^{-1}
(\partial_\tau r)\sin\phi$. Equation \eqref{eq:cap-meridian-components}
 give
\begin{equation}\label{eq:cap-moment-normal-derivative}
\begin{aligned}
 -\frac{r^{n-2}\sin^{d-2}\!\phi}{\partial_\tau r}\partial_\nu Z_e
 &=r^{n-2}\sin^{d-2}\phi
   \left(\alpha_{\mathrm{FS}}r\phi'\cos\phi
         +\alpha_{\mathrm{FS}}^{-1}\sin\phi\right)\\
 &=\left(\frac{\alpha_{\mathrm{FS}}}{d-1}
          r^{n-1}\sin^{d-1}\phi\right)'.
\end{aligned}
\end{equation}
Equations \eqref{eq:cap-moment-product} and
\eqref{eq:cap-moment-normal-derivative} give
\begin{equation}\label{eq:cap-moment-derivative}
\frac{d}{dr}\left(
 r^{n-2}\sin^{d-2}\!\phi\,\partial_\tau Z_e
 -\frac{\alpha_{\mathrm{FS}}\lambda_E}{d-1}r^{n-1}\sin^{d-1}\phi
 \right)=\frac{n-d}{d-1}\,r^{n-3}\sin^{d-2}\!\phi\,
          (\partial_\tau r)\cos\phi .
\end{equation}
Since $d\ge3$ and $|\partial_\tau Z_e|\le\alpha_{\mathrm{FS}}^{-1}$,
both boundary terms vanish at a positive endpoint by
\eqref{eq:cap-meridian-endpoints}. If $r_-=0$, the factors $r^{n-2}$
and $r^{n-1}$ also vanish. Since
$r^{n-3}\in L^1(0,r_+)$, integration of \eqref{eq:cap-moment-derivative} gives
\[
 0=\frac{n-d}{d-1}
   \int_{r_-}^{r_+}r^{n-3}\sin^{d-2}\!\phi\,(\partial_\tau r)\cos\phi\,dr.
\]
Since $n>d$, this gives \eqref{eq:cap-profile-moment}. If
$\cos\phi\ge0$ throughout $(r_-,r_+)$, positivity of
$r^{n-3}\sin^{d-2}\!\phi\,\partial_\tau r$ and continuity would force
$\phi\equiv\pi/2$, contrary to \eqref{eq:cap-meridian-endpoints} at $r_+>0$.
\end{proof}

\begin{lemma}\label{lem:cap-inner-reflection}
In the nonball case of Lemma~\ref{lem:cap-meridian}, suppose that
$\phi(r_0)>\pi/2$ for some $r_0\in(r_-,r_+)$. Then
\begin{equation}\label{eq:inner-cap-hemisphere}
 |E_r|\ge\frac{\omega_{d-1}}2
 \quad\text{for a.e. }r\in(0,r_0),\qquad
 \phi(r)\ge\frac{\pi}{2}
 \quad\text{for every }r\in(r_-,r_0).
\end{equation}
\end{lemma}

\begin{proof}
Let $\mathcal Rz:=z-2(z\cdot e)e$ be reflection across $e^\perp$, and set
\[
 E':=E\cup\bigl(B_{r_0}(0)\cap\mathcal R(\Omega\setminus E)\bigr).
\]
For a.e.\ $0<r<r_0$, the slices of $E$ and $\mathcal R(\Omega\setminus E)$
are nested caps up to null sets, of measures $|E_r|$ and
$\omega_{d-1}-|E_r|$. Thus
\begin{equation}\label{eq:cap-reflected-slice-gain}
 |(E')_r|-|E_r|
 =\max\{|E_r|,\omega_{d-1}-|E_r|\}-|E_r|
 =(\omega_{d-1}-2|E_r|)_+.
\end{equation}
In $B_{r_0}(0)$, $E'$ agrees with $E\cup\mathcal R(\Omega\setminus E)$
and $\mathcal R(\Omega\setminus E')$ with $E\cap\mathcal R(\Omega\setminus E)$,
since $\mathcal R$ preserves the ball and $\mathcal R^2$ is the identity.
The local union--intersection inequality \cite[Lemma~12.22]{M12} with
the radial weight in \eqref{eq:cone-pullback-volume-perimeter} gives
\begin{equation}\label{eq:cap-reflected-perimeter}
\begin{aligned}
 2P_{\alpha_{\mathrm{FS}}}(E';B_{r_0}(0))
 &=P_{\alpha_{\mathrm{FS}}}(E';B_{r_0}(0))
   +P_{\alpha_{\mathrm{FS}}}\bigl(E\cap\mathcal R(\Omega\setminus E);B_{r_0}(0)\bigr)\\
 &\le P_{\alpha_{\mathrm{FS}}}(E;B_{r_0}(0))
   +P_{\alpha_{\mathrm{FS}}}\bigl(\mathcal R(\Omega\setminus E);B_{r_0}(0)\bigr)\\
 &=2P_{\alpha_{\mathrm{FS}}}(E;B_{r_0}(0)).
\end{aligned}
\end{equation}
By continuity, $\phi>\pi/2$ for $|r-r_0|<\delta$, with
$(r_0-\delta,r_0+\delta)\subset(r_-,r_+)$. On this interval
\eqref{eq:cap-meridian-slices} gives
\begin{equation}\label{eq:cap-reflected-containment}
 \bigl(\mathcal R(\Omega\setminus E^{(1)})\bigr)_r
 =\{\theta:\vartheta(\theta)\le\pi-\phi(r)\}
 \subset\{\theta:\vartheta(\theta)<\phi(r)\}=(E^{(1)})_r.
\end{equation}
Since $\mathcal R$ preserves $\mu$ and $\mu(E\triangle E^{(1)})=0$,
\eqref{eq:cap-reflected-containment} gives $E'=E$ a.e.\ near
$\partial B_{r_0}(0)$. The sets also agree outside the ball, so their
perimeter measures coincide on $\Omega\setminus B_{r_0}(0)$. By
\eqref{eq:cap-reflected-perimeter} and \eqref{eq:cap-reflected-slice-gain},
\begin{align}
 P_{\alpha_{\mathrm{FS}}}(E';\Omega)-P_{\alpha_{\mathrm{FS}}}(E;\Omega)
 &=P_{\alpha_{\mathrm{FS}}}(E';B_{r_0}(0))-P_{\alpha_{\mathrm{FS}}}(E;B_{r_0}(0))\le0,
 \label{eq:cap-reflected-competitor-perimeter}\\
 \mu(E')-\mu(E)
 &=\int_0^{r_0}r^{n-1}(\omega_{d-1}-2|E_r|)_+\,dr\ge0.
 \label{eq:cap-reflected-competitor-volume}
\end{align}
For $s:=(\mu(E)/\mu(E'))^{1/n}\le1$, \eqref{eq:cap-dilation-laws}
gives $\mu(sE')=\mu(E)$. Minimality and \eqref{eq:cap-reflected-competitor-perimeter}
give
\[
 0<P_{\alpha_{\mathrm{FS}}}(E;\Omega)
 \le P_{\alpha_{\mathrm{FS}}}(sE';\Omega)
 =s^{n-1}P_{\alpha_{\mathrm{FS}}}(E';\Omega)
 \le s^{n-1}P_{\alpha_{\mathrm{FS}}}(E;\Omega).
\]
Hence $s=1$ and \eqref{eq:cap-reflected-competitor-volume} gives
$|E_r|\ge\omega_{d-1}/2$ for a.e.\ $r<r_0$. On $(r_-,r_0)$,
\eqref{eq:cap-opening-angle} and \eqref{eq:cap-meridian-slices} give
\[
 |(E^{(1)})_r|-\frac{\omega_{d-1}}2
 =|\mathbb S^{d-2}|\int_{\pi/2}^{\phi(r)}\sin^{d-2}\!t\,dt.
\]
The left side is nonnegative a.e.\ by \eqref{eq:cap-representative-slices},
and the right side is continuous. Since $\sin^{d-2}t>0$ on $(0,\pi)$,
$\phi(r)\ge\pi/2$ throughout $(r_-,r_0)$.
\end{proof}

\begin{lemma}\label{lem:cap-interior-ball}
In the nonball case of Lemma~\ref{lem:cap-meridian}, $r_->0$,
$B_{r_-}(0)\subset E^{(1)}$, and the meridian has endpoints
$p_-=-r_-e$ and $p_+=r_+e$ with $\psi(p_\pm)=r_\pm>0$.
\end{lemma}

\begin{proof}
By Lemmas~\ref{lem:cap-angular-moment} and~\ref{lem:cap-inner-reflection},
there is $r_0\in(r_-,r_+)$ with $\phi(r_0)>\pi/2$ and
$\phi\ge\pi/2$ on $(r_-,r_0)$.

Suppose $r_-=0$. Since $\cot\phi\le0$ on $(0,r_0]$,
\eqref{eq:cap-slope-ode} of Lemma~\ref{lem:cone-meridian-system} gives
\begin{equation}\label{eq:cap-inner-flux}
 (r^{n-2}\psi)'
 =\lambda_E r^{n-1}
  -\frac{d-2}{\alpha_{\mathrm{FS}}}r^{n-3}\cot\phi\,\sqrt{r^2-\psi^2}
 \ge\lambda_E r^{n-1}.
\end{equation}
As $|\psi|\le r$, $r^{n-2}\psi(r)\to0$ at the origin. Integrating
\eqref{eq:cap-inner-flux} gives $\psi(r)\ge\lambda_E r^2/n>0$. Thus \eqref{eq:cap-angle-ode} gives
\begin{equation}\label{eq:cap-inner-angle-sign}
 \phi'(r)<0,\qquad
 \phi(r)\ge\phi(r_0)>\frac{\pi}{2},\qquad
 \cot\phi(r)\le\cot\phi(r_0)<0
 \qquad(0<r\le r_0).
\end{equation}
Put $c_{\phi}:=-(d-2)\cot\phi(r_0)/\alpha_{\mathrm{FS}}>0$.
By \eqref{eq:cap-meridian-speed}, $0<\psi<r$. Using
$\sqrt{r^2-\psi^2}\ge r-\psi$ and \eqref{eq:cap-inner-angle-sign}
in \eqref{eq:cap-slope-ode} gives
\[
 r\psi'+(n-2)\psi
 =\lambda_E r^2-\frac{d-2}{\alpha_{\mathrm{FS}}}
   \cot\phi\,\sqrt{r^2-\psi^2}
 \ge c_{\phi}(r-\psi).
\]
Thus $(r^{n-2+c_{\phi}}\psi)'\ge c_{\phi}r^{n-2+c_{\phi}}$.
Integrating from $0$ to $r$ and using $|\psi|\le r$, we obtain
\begin{equation}\label{eq:cap-inner-slope-bound}
 \frac{\psi(r)}{r}\ge\frac{c_{\phi}}{n-1+c_{\phi}}
 \qquad(0<r\le r_0).
\end{equation}
By \eqref{eq:cap-angle-ode} and \eqref{eq:cap-inner-slope-bound},
\[
 \phi(r)-\phi(r_0)
 =\int_r^{r_0}\frac{\psi(s)}
    {\alpha_{\mathrm{FS}}s\sqrt{s^2-\psi(s)^2}}\,ds
 \ge \frac{c_{\phi}}{\alpha_{\mathrm{FS}}(n-1+c_{\phi})}
     \log\frac{r_0}{r}
 \to+\infty
 \qquad(r\downarrow0).
\]
Since $\phi<\pi$, we conclude that $r_->0$.

The punctured ball $B_{r_-}(0)\setminus\{0\}$ is connected and disjoint
from $\partial E^{(1)}$. By \eqref{eq:inner-cap-hemisphere},
\[
 \mu(E^{(1)}\cap B_{r_-}(0))
 =\int_0^{r_-}r^{n-1}|E_r|\,dr
 \ge\tfrac12\mu(B_{r_-}(0))>0.
\]
Thus the punctured ball lies in $\operatorname{int}E^{(1)}$, and
$0\in E^{(1)}$ by \eqref{eq:cap-density-one}. The connected exterior
$\{r>r_+\}$ avoids $\partial E^{(1)}$ and contains points outside the
bounded set $E^{(1)}$. Hence
\begin{equation}\label{eq:cap-inner-outer-regions}
 B_{r_-}(0)\subset E^{(1)},\qquad
 \{r>r_+\}\cap\overline{E^{(1)}}=\varnothing.
\end{equation}
By \eqref{eq:cap-meridian-slices} and $0<\phi<\pi$,
\begin{equation}\label{eq:cap-axial-membership}
 re\in\operatorname{int}E^{(1)},\qquad -re\notin\overline{E^{(1)}}
 \qquad(r_-<r<r_+).
\end{equation}
Equations \eqref{eq:cap-inner-outer-regions} and
\eqref{eq:cap-axial-membership} put $-r_-e,r_+e$ on $\Sigma$.
By \eqref{eq:cap-meridian-endpoints},
\[
 p_-=-r_-e,\qquad p_+=r_+e,\qquad
 \lim_{r\downarrow r_-}\phi(r)=\pi,\qquad
 \lim_{r\uparrow r_+}\phi(r)=0.
\]
On both rays, increasing $r$ crosses from $\operatorname{int}E^{(1)}$
to its exterior by \eqref{eq:cap-inner-outer-regions} and
\eqref{eq:cap-axial-membership}. Axial symmetry gives a radial normal,
so $\nu(p_\pm)=\partial_r|_{p_\pm}$ and $\psi(p_\pm)=r_\pm>0$
by \eqref{eq:cap-dilation-speed}.
\end{proof}

\begin{proof}[Proof of Proposition~\ref{prop:cap-strict-starshaped}]
For $E^{(1)}=B_{(n-1)/\lambda_E}(0)$, $\nu=\partial_r$ and
$\psi=R(\theta)=(n-1)/\lambda_E>0$ by \eqref{eq:cap-dilation-speed}.

Otherwise $p_\pm\in\Omega_+$ by Lemma~\ref{lem:cap-interior-ball}.
By Lemma~\ref{lem:cap-positive-connected}, $r(\Omega_+)$ is connected.
Hence
\begin{equation}\label{eq:cap-positive-radii}
 \{r_-,r_+\}\subset r(\Omega_+)\subset[r_-,r_+],
 \qquad r(\Omega_+)=[r_-,r_+].
\end{equation}
For each $r\in(r_-,r_+)$, choose $z_r\in\Omega_+\cap\partial B_r(0)$
by \eqref{eq:cap-positive-radii}. The latitude in \eqref{eq:cap-meridian-slices}
is a single $G_e$-orbit. This group, defined in \eqref{eq:cap-polar-angle},
preserves $E^{(1)}$, $g_{\alpha_{\mathrm{FS}}}$, and $r\partial_r$,
so \eqref{eq:cap-dilation-speed} gives
\[
 \psi(Oz_r)=\psi(z_r)>0\qquad(O\in G_e).
\]
By \eqref{eq:cap-meridian-endpoints}, only $p_\pm$ remain, and
$\psi(p_\pm)>0$. Thus $\psi>0$ on $\Sigma$. For every local
defining function $F$, \eqref{eq:cap-defining-function-speed} gives
\begin{equation}\label{eq:cap-radial-transversality}
 \partial_r F
 =\frac{|\nabla_{\alpha_{\mathrm{FS}}}F|}{r}
   \,\psi>0
 \qquad\text{on }\Sigma\cap U.
\end{equation}
Fix $\theta\in\mathbb S^{d-1}$, choose $0<\varepsilon<r_-<r_+<M$,
and set $S_\theta:=\{s>0:s\theta\in\Sigma\}$. For each $s\in S_\theta$,
\eqref{eq:cap-radial-transversality} and $E^{(1)}\cap U=\{F<0\}$ give
\begin{equation}\label{eq:cap-ray-exit}
 \mathbf1_{E^{(1)}}((s-h)\theta)=1,\qquad
 \mathbf1_{E^{(1)}}((s+h)\theta)=0
 \qquad\text{for all sufficiently small }h>0.
\end{equation}
By \eqref{eq:cap-radial-transversality}, $S_\theta$ is discrete.
It is closed in $[\varepsilon,M]$, hence finite. The indicator
is constant between crossings and drops by one at each by
\eqref{eq:cap-ray-exit}. Summing and using \eqref{eq:cap-inner-outer-regions}
gives
\[
 \operatorname{card}(S_\theta)
 =\mathbf1_{E^{(1)}}(\varepsilon\theta)-\mathbf1_{E^{(1)}}(M\theta)=1.
\]
Write $R(\theta)$ for this unique radius. The implicit function
theorem and \eqref{eq:cap-radial-transversality} give smoothness in
each angular chart. Since $B_{r_-}(0)\subset E^{(1)}$,
\[
 R\in C^\infty(\mathbb S^{d-1}),\qquad R\ge r_->0,\qquad
 \partial E^{(1)}=\{R(\theta)\theta:\theta\in\mathbb S^{d-1}\}.
\]
This completes the proof.
\end{proof}

\section{Rigidity of stationary radial graphs}
\label{sec:graph-rigidity}

Fix $n>d\ge3$. All derivatives ($\nabla_\theta$, $\nabla_\theta^2$,
$\Delta_\theta$, $\operatorname{div}_\theta$) and inner products are taken
in the round metric on $\mathbb S^{d-1}$, with normalized surface measure $d\sigma$.

\begin{proposition}\label{prop:stationary-graphs}
Let $n>d\ge3$. If a smooth positive radial graph $r=R(\theta)$ on
$\mathbb S^{d-1}$ is stationary for $P_{\alpha_{\mathrm{FS}}}(\cdot;\Omega)$ at fixed
weighted volume, then $R$ is constant.
\end{proposition}

\subsection{The stationary equation and a volume-preserving variation}

For $h\in C^\infty(\mathbb S^{d-1})$, write $R=e^h$ and let
$E[e^h]$ be the subgraph in \eqref{eq:radial-subgraph}. Set
\begin{equation}\label{eq:graph-slope}
 m_h:=\bigl(1+\alpha_{\mathrm{FS}}^{-2}|\nabla_\theta h|^2\bigr)^{1/2}.
\end{equation}
By \eqref{eq:radial-graph-perimeter} and radial integration in
\eqref{eq:cone-metric-measure},
\begin{align}
 \mathcal A(h)
 &:=\frac{P_{\alpha_{\mathrm{FS}}}(E[e^h];\Omega)}{\omega_{d-1}}
   =\int_{\mathbb S^{d-1}}e^{(n-1)h}m_h\,d\sigma,
 \label{eq:def-Pv}\\
 \mathcal V(h)
 &:=\frac{\mu(E[e^h])}{\omega_{d-1}}
   =\frac1n\int_{\mathbb S^{d-1}}e^{nh}\,d\sigma.
 \label{eq:def-Vv}
\end{align}

At a stationary graph, $\mathcal A'(h)$ is proportional to $\mathcal V'(h)$.
The direction $u=1$ fixes the multiplier, since $\mathcal A'(h)[1]=(n-1)\mathcal A(h)$
and $\mathcal V'(h)[1]=n\mathcal V(h)>0$. Thus
\begin{equation}\label{eq:graph-constrained-stationarity}
 \mathcal A'(h)[u]=\frac{(n-1)\mathcal A(h)}{n\mathcal V(h)}\,\mathcal V'(h)[u]
 \qquad\text{for every }u\in C^\infty(\mathbb S^{d-1}).
\end{equation}
By \eqref{eq:def-Pv} and \eqref{eq:def-Vv},
\eqref{eq:graph-constrained-stationarity} reads
\begin{equation}\label{eq:graph-weak-stationarity}
 \int_{\mathbb S^{d-1}}e^{(n-1)h}
 \left((n-1)m_hu
       +\alpha_{\mathrm{FS}}^{-2}
        \frac{\nabla_\theta h\cdot\nabla_\theta u}{m_h}\right)d\sigma
 =\frac{(n-1)\mathcal A(h)}{n\mathcal V(h)}\int_{\mathbb S^{d-1}}e^{nh}u\,d\sigma.
\end{equation}
Integrating \eqref{eq:graph-weak-stationarity} by parts gives
\begin{equation}\label{eq:graph-stationarity-divergence}
 (n-1)e^{(n-1)h}m_h
 -\alpha_{\mathrm{FS}}^{-2}\operatorname{div}_{\theta}\!\left(
   e^{(n-1)h}\frac{\nabla_\theta h}{m_h}\right)
 =\frac{(n-1)\mathcal A(h)}{n\mathcal V(h)}e^{nh}.
\end{equation}

Define
\begin{equation}\label{eq:graph-T}
 T:=m_h\operatorname{div}_{\theta}\!\left(\frac{\nabla_\theta h}{m_h}\right)
   =\Delta_\theta h-\frac{\nabla_\theta^2h(\nabla_\theta h,\nabla_\theta h)}
                    {\alpha_{\mathrm{FS}}^2+|\nabla_\theta h|^2}.
\end{equation}
The product rule and \eqref{eq:graph-T} give
\begin{equation}\label{eq:graph-weighted-divergence}
 \operatorname{div}_{\theta}\!\left(e^{(n-1)h}\frac{\nabla_\theta h}{m_h}\right)
 =\frac{e^{(n-1)h}}{m_h}\bigl(T+(n-1)|\nabla_\theta h|^2\bigr).
\end{equation}
By \eqref{eq:graph-stationarity-divergence}, \eqref{eq:graph-weighted-divergence},
and $m_h^2-\alpha_{\mathrm{FS}}^{-2}|\nabla_\theta h|^2=1$,
\[
 (n-1)-\alpha_{\mathrm{FS}}^{-2}T=\frac{(n-1)\mathcal A(h)}{n\mathcal V(h)}e^h m_h.
\]

Take
\[
 \zeta:=e^{-nh}\operatorname{div}_{\theta}\!\left(
          e^{(n-1)h}\frac{\nabla_\theta h}{m_h}\right).
\]
By \eqref{eq:graph-weighted-divergence} and the divergence theorem,
\begin{equation}\label{eq:graph-volume-test}
 \zeta=\frac{e^{-h}}{m_h}\bigl(T+(n-1)|\nabla_\theta h|^2\bigr),
 \qquad
 \mathcal V'(h)[\zeta]
 =\int_{\mathbb S^{d-1}}e^{nh}\zeta\,d\sigma=0.
\end{equation}
The path
\begin{equation}\label{eq:graph-volume-preserving-path}
 h_\varepsilon:=h+\varepsilon\zeta-\frac1n
       \log\!\left(\frac{\mathcal V(h+\varepsilon\zeta)}{\mathcal V(h)}\right)
\end{equation}
preserves volume, since $\mathcal V(k+c)=e^{nc}\mathcal V(k)$ for
constant $c$. By \eqref{eq:graph-volume-test},
\[
 \mathcal V(h_\varepsilon)=\mathcal V(h),\qquad
 \left.\frac{d}{d\varepsilon}\right|_{\varepsilon=0}h_\varepsilon
 =\zeta-\frac{\mathcal V'(h)[\zeta]}{n\mathcal V(h)}
 =\zeta.
\]

Integration by parts in the first variation of \eqref{eq:def-Pv},
using \eqref{eq:graph-weighted-divergence}, gives
\begin{equation}\label{eq:graph-area-first-variation}
 \mathcal A'(h)[u]
 =\int_{\mathbb S^{d-1}}\frac{e^{(n-1)h}}{m_h}
   \bigl((n-1)-\alpha_{\mathrm{FS}}^{-2}T\bigr)u\,d\sigma.
\end{equation}
Set
\begin{equation}\label{eq:graph-angular-weight}
 \mathcal W_h:=\frac{e^{(n-2)h}}{m_h^2}.
\end{equation}
Stationarity along \eqref{eq:graph-volume-preserving-path} gives
$\mathcal A'(h)[\zeta]=0$. By \eqref{eq:graph-volume-test},
\eqref{eq:graph-area-first-variation}, \eqref{eq:graph-angular-weight},
and $\alpha_{\mathrm{FS}}^2(n-1)=d-1$,
\begin{equation}\label{eq:angular-stationary-test}
 0=\mathcal A'(h)[\zeta]
 =-\alpha_{\mathrm{FS}}^{-2}\int_{\mathbb S^{d-1}}
   \mathcal W_h\bigl(T-(d-1)\bigr)
   \bigl(T+(n-1)|\nabla_\theta h|^2\bigr)\,d\sigma.
\end{equation}

\subsection{The nonlinear angular inequality and graph rigidity}

In a local orthonormal frame $(e_i)_{i=1}^{d-1}$, write
$h_i=\langle\nabla_\theta h,e_i\rangle$ and $h_{ij}=\nabla_\theta^2h(e_i,e_j)$.
Identifying the Hessian with its self-adjoint endomorphism, we have
\[
 \|\nabla_\theta^2h\|_{\mathrm{HS}}^2=\sum_{i,j}h_{ij}^2,\qquad
 |\nabla_\theta^2h\,\nabla_\theta h|^2=\sum_i\left(\sum_j h_{ij}h_j\right)^2,
\]
where $\|\cdot\|_{\mathrm{HS}}$ denotes the Hilbert--Schmidt norm.
Define
\begin{equation}\label{eq:graph-angular-hessian}
 S:=\|\nabla_\theta^2h\|_{\mathrm{HS}}^2
   -\frac{2|\nabla_\theta^2h\,\nabla_\theta h|^2}{\alpha_{\mathrm{FS}}^2+|\nabla_\theta h|^2}
   +\frac{\nabla_\theta^2h(\nabla_\theta h,\nabla_\theta h)^2}
          {(\alpha_{\mathrm{FS}}^2+|\nabla_\theta h|^2)^2}.
\end{equation}

Set
\begin{equation}\label{eq:graph-P}
 P_h e_i:=e_i-\frac{h_i\nabla_\theta h}
              {\alpha_{\mathrm{FS}}^2m_h(m_h+1)}.
\end{equation}
Squaring \eqref{eq:graph-P} and using
$|\nabla_\theta h|^2=\alpha_{\mathrm{FS}}^2(m_h^2-1)$ gives
\begin{equation}\label{eq:graph-P-square}
 P_h^2 e_i=e_i-\frac{h_i\nabla_\theta h}
              {\alpha_{\mathrm{FS}}^2+|\nabla_\theta h|^2}.
\end{equation}
Cyclicity of the trace and \eqref{eq:graph-P-square} give
\[
\begin{aligned}
 \operatorname{tr}(P_h\nabla_\theta^2h\,P_h)
 &=\operatorname{tr}(\nabla_\theta^2h\,P_h^2)=T,\qquad
 \|P_h\nabla_\theta^2h\,P_h\|_{\mathrm{HS}}^2
 &=\operatorname{tr}(\nabla_\theta^2h\,P_h^2\nabla_\theta^2h\,P_h^2)=S,
\end{aligned}
\]
where $T$ and $S$ are given by \eqref{eq:graph-T} and \eqref{eq:graph-angular-hessian}.
 Cauchy--Schwarz applied to the
eigenvalues of $P_h\nabla_\theta^2h\,P_h$ gives $T^2\le(d-1)S$.

\begin{proposition}\label{pro:angular-bellman}
Let $n>d\ge3$. For every $h\in C^\infty(\mathbb S^{d-1})$, with $T$, $\mathcal W_h$, and $S$ as in \eqref{eq:graph-T},
\eqref{eq:graph-angular-weight}, and \eqref{eq:graph-angular-hessian},
\begin{equation}\label{eq:angular-bellman}
 \mathcal E_{\mathrm{ang}}[h]:=\int_{\mathbb S^{d-1}}\mathcal W_h
 \left(S-\frac{T^2}{n-1}-\frac{n-d}{n-1}|\nabla_\theta h|^2\right)d\sigma\ge0,
\end{equation}
with equality if and only if $h$ is constant.
\end{proposition}

\begin{lemma}\label{lem:angular-divergences}
For $h\in C^\infty(\mathbb S^{d-1})$, set $J_h:=(\Delta_\theta h)\nabla_\theta h-\nabla_\theta^2h\,\nabla_\theta h$.
Define
\begin{equation}\label{eq:graph-hessian-contractions}
 \begin{aligned}
 D_h&:=|\nabla_\theta h|^2\Delta_\theta h-\nabla_\theta^2h(\nabla_\theta h,\nabla_\theta h),\\
 B_h&:=\nabla_\theta^2h(\nabla_\theta h,\nabla_\theta h)\Delta_\theta h-|\nabla_\theta^2h\,\nabla_\theta h|^2.
\end{aligned}
\end{equation}
For every $\mathcal M\in C^1([0,\infty))$, evaluate $\mathcal M$ and $\mathcal M'$
at $|\nabla_\theta h|^2$. Then
\begin{align}
 0&=\int_{\mathbb S^{d-1}}\mathcal W_h
 \Bigl\{\mathcal M\bigl[T^2-S-(d-2)|\nabla_\theta h|^2\bigr]
       +(n-2)\mathcal M D_h+2\mathcal M'B_h\Bigr\}\,d\sigma,
 \label{eq:bellman-zero-one}\\
 0&=\int_{\mathbb S^{d-1}}e^{(n-2)h}
 \Bigl\{\mathcal M\Delta_\theta h+2\mathcal M'\nabla_\theta^2h(\nabla_\theta h,\nabla_\theta h)
       +(n-2)\mathcal M|\nabla_\theta h|^2\Bigr\}\,d\sigma.
 \label{eq:bellman-zero-two}
\end{align}
\end{lemma}

\begin{proof}
Since $\nabla_\theta(|\nabla_\theta h|^2)=2\nabla_\theta^2 h\,\nabla_\theta h$,
\eqref{eq:graph-hessian-contractions} gives
\begin{equation}\label{eq:graph-J-contractions}
 J_h\cdot\nabla_\theta h=D_h,\quad
 J_h\cdot\nabla_\theta\bigl(|\nabla_\theta h|^2\bigr)=2B_h,\quad
 \nabla_\theta\bigl(|\nabla_\theta h|^2\bigr)\cdot\nabla_\theta h=2\nabla_\theta^2h(\nabla_\theta h,\nabla_\theta h).
\end{equation}
The spherical Bochner formula \cite[proof of Lemma~5.2]{DEL16} gives
\begin{equation}\label{eq:graph-bochner-divergence}
\begin{aligned}
 \operatorname{div}_{\theta}J_h
 &=(\Delta_\theta h)^2+\nabla_\theta h\cdot\nabla_\theta\Delta_\theta h-\tfrac12\Delta_\theta\bigl(|\nabla_\theta h|^2\bigr)\\
 &=(\Delta_\theta h)^2-\|\nabla_\theta^2h\|_{\mathrm{HS}}^2-(d-2)|\nabla_\theta h|^2.
\end{aligned}
\end{equation}

Equations \eqref{eq:graph-T},
\eqref{eq:graph-angular-hessian}, and \eqref{eq:graph-hessian-contractions}
give
\begin{equation}\label{eq:graph-trace-square-difference}
 T^2-S
 =(\Delta_\theta h)^2-\|\nabla_\theta^2h\|_{\mathrm{HS}}^2
  -\frac{2B_h}{\alpha_{\mathrm{FS}}^2+|\nabla_\theta h|^2}.
\end{equation}
Differentiating $m_h^{-2}$ and using \eqref{eq:graph-J-contractions},
\eqref{eq:graph-bochner-divergence}, and \eqref{eq:graph-trace-square-difference},
we obtain
\begin{equation}\label{eq:graph-weighted-bochner}
 \operatorname{div}_{\theta}(m_h^{-2}J_h)
 =m_h^{-2}\left(\operatorname{div}_{\theta}J_h
                 -\frac{2B_h}{\alpha_{\mathrm{FS}}^2+|\nabla_\theta h|^2}\right)=m_h^{-2}\bigl(T^2-S-(d-2)|\nabla_\theta h|^2\bigr).
\end{equation}

By \eqref{eq:graph-angular-weight}, \eqref{eq:graph-weighted-bochner},
and \eqref{eq:graph-J-contractions},
\begin{equation}\label{eq:graph-first-multiplier-divergence}
\begin{aligned}
 \operatorname{div}_{\theta}(\mathcal W_h\mathcal M J_h)
 &=e^{(n-2)h}\mathcal M\,\operatorname{div}_{\theta}(m_h^{-2}J_h)
   +\mathcal W_h J_h\cdot
    \bigl((n-2)\mathcal M\nabla_\theta h+\mathcal M'\nabla_\theta\bigl(|\nabla_\theta h|^2\bigr)\bigr)\\
 &=\mathcal W_h
   \Bigl\{\mathcal M\bigl[T^2-S-(d-2)|\nabla_\theta h|^2\bigr]
         +(n-2)\mathcal M D_h+2\mathcal M'B_h\Bigr\}.
\end{aligned}
\end{equation}
Likewise,
\begin{equation}\label{eq:graph-second-multiplier-divergence}
\begin{aligned}
 \operatorname{div}_{\theta}(e^{(n-2)h}\mathcal M\nabla_\theta h)
 &=e^{(n-2)h}\bigl(\mathcal M\Delta_\theta h+\mathcal M'\nabla_\theta\bigl(|\nabla_\theta h|^2\bigr)\cdot\nabla_\theta h+(n-2)\mathcal M |\nabla_\theta h|^2\bigr)\\
 &=e^{(n-2)h}
   \bigl(\mathcal M\Delta_\theta h+2\mathcal M'\nabla_\theta^2h(\nabla_\theta h,\nabla_\theta h)+(n-2)\mathcal M |\nabla_\theta h|^2\bigr).
\end{aligned}
\end{equation}
Integrating \eqref{eq:graph-first-multiplier-divergence} and
\eqref{eq:graph-second-multiplier-divergence} gives
\eqref{eq:bellman-zero-one} and \eqref{eq:bellman-zero-two}.
\end{proof}

\begin{proof}[Proof of Proposition~\ref{prop:stationary-graphs}]
Let $h=\log R$. All integrals are over $\mathbb S^{d-1}$.
By \eqref{eq:graph-slope}, \eqref{eq:graph-T}, and \eqref{eq:graph-hessian-contractions},
\begin{equation}\label{eq:graph-G-identities}
 |\nabla_\theta h|^2T-D_h
 =\frac{\nabla_\theta^2h(\nabla_\theta h,\nabla_\theta h)}{m_h^2}
 =\alpha_{\mathrm{FS}}^2(\Delta_\theta h-T).
\end{equation}
Take $\mathcal M=1$ in \eqref{eq:bellman-zero-one} of
Lemma~\ref{lem:angular-divergences}. Then
\begin{equation}\label{eq:graph-integrated-hessian}
 \int\mathcal W_h S\,d\sigma
 =\int\mathcal W_h
   \bigl[T^2-(d-2)|\nabla_\theta h|^2+(n-2)D_h\bigr]\,d\sigma.
\end{equation}
In \eqref{eq:bellman-zero-two}, take $\mathcal M(s)=(1+\alpha_{\mathrm{FS}}^{-2}s)^{-1}$.
At $s=|\nabla_\theta h|^2$, $\mathcal M=m_h^{-2}$ and
$\mathcal M'=-\alpha_{\mathrm{FS}}^{-2}m_h^{-4}$. By \eqref{eq:graph-G-identities},
\[
 m_h^{-2}\Delta_\theta h-2\alpha_{\mathrm{FS}}^{-2}m_h^{-4}
 \nabla_\theta^2h(\nabla_\theta h,\nabla_\theta h)
 =m_h^{-2}\bigl[T-\alpha_{\mathrm{FS}}^{-2}(|\nabla_\theta h|^2T-D_h)\bigr].
\]
Thus \eqref{eq:bellman-zero-two}, \eqref{eq:graph-angular-weight}, and
$(n-1)\alpha_{\mathrm{FS}}^2=d-1$ give
\begin{equation}\label{eq:graph-integrated-G}
 (n-1)\int\mathcal W_h D_h\,d\sigma
 =\int\mathcal W_h\bigl[(n-1)|\nabla_\theta h|^2T-(d-1)(T+(n-2)|\nabla_\theta h|^2)\bigr]\,d\sigma .
\end{equation}
By \eqref{eq:angular-bellman}, \eqref{eq:graph-integrated-hessian},
and \eqref{eq:graph-integrated-G},
\begin{equation}\label{eq:graph-angular-factorization}
\begin{aligned}
 \frac{n-1}{n-2}\mathcal E_{\mathrm{ang}}[h]
 &=\int\mathcal W_h
   \bigl[T^2+(n-1)D_h-(d-1)|\nabla_\theta h|^2\bigr]\,d\sigma\\
 &=\int\mathcal W_h
   \bigl[T^2+(n-1)|\nabla_\theta h|^2T-(d-1)T-(d-1)(n-1)|\nabla_\theta h|^2\bigr]\,d\sigma\\
 &=\int\mathcal W_h
   \bigl(T-(d-1)\bigr)\bigl(T+(n-1)|\nabla_\theta h|^2\bigr)\,d\sigma.
\end{aligned}
\end{equation}

Since $\mathcal V'(h)[\zeta]=0$, stationarity gives
$\mathcal A'(h)[\zeta]=0$ by \eqref{eq:graph-constrained-stationarity}.
By \eqref{eq:angular-stationary-test} and \eqref{eq:graph-angular-factorization},
\[
 0=\mathcal A'(h)[\zeta]
  =-\frac{n-1}{\alpha_{\mathrm{FS}}^2(n-2)}\,\mathcal E_{\mathrm{ang}}[h].
\]
Equality in Proposition~\ref{pro:angular-bellman} forces $h$ and hence $R=e^h$ to be constant.
\end{proof}

\subsection{Proof of Proposition~\ref{pro:angular-bellman}}
At $\alpha=\alpha_{\mathrm{FS}}$, the angular part of the weighted
cone's $n$-Bakry--\'Emery Ricci tensor, pulled back to the unit
sphere, is
\[
 \bigl[(d-2)-(n-2)\alpha_{\mathrm{FS}}^2\bigr]g_{\mathbb S^{d-1}}
 =-\frac{n-d}{n-1}g_{\mathbb S^{d-1}}.
\]
The corresponding term in \eqref{eq:angular-bellman} is
$-\frac{n-d}{n-1}|\nabla_\theta h|^2$. To compensate for it, we add
\eqref{eq:bellman-zero-one} and \eqref{eq:bellman-zero-two}, in the
spirit of the carr\'e du champ argument in
\cite[proof of Lemma~5.2]{DEL16} for $p=2$.

As a function of the multiplier argument $s\ge0$ in
\eqref{eq:bellman-multipliers}, the part independent of the Hessian is
\[
 -\frac{n-d}{n-1}\frac{\alpha_{\mathrm{FS}}^{-2}s^2}
                         {\tfrac89+8(1-\kappa)^2+\alpha_{\mathrm{FS}}^{-2}s}
 +(n-2)\gamma\frac{s^2}{\sqrt{1+\alpha_{\mathrm{FS}}^{-2}s}}.
\]
The value $\mathcal M_1(0)=-(n-d)/[(n-1)(d-2)]$ cancels the linear
term at zero, and $s(1+\alpha_{\mathrm{FS}}^{-2}s)\mathcal M_1'(s)$
stays bounded as $s\to\infty$. The positive term from $\mathcal M_2$
grows like $s^{3/2}$ and dominates the negative term of order $s$ at
infinity. Completing squares in the Hessian reduces the case $s>0$
to Lemma~\ref{lem:angular-sign}.

\begin{proof}[Proof of Proposition~\ref{pro:angular-bellman}]
Fix $h\in C^\infty(\mathbb S^{d-1})$, with $m_h$ as in
\eqref{eq:graph-slope}. Set
\begin{equation}\label{eq:bellman-parameters}
 \kappa:=\frac{n-d}{n-2},\qquad
 \beta:=\frac{n-d}{(n-1)(d-2)},\qquad
 \gamma:=\frac{(n-2)(n-d)}{(n-1)(d+1)},
\end{equation}
and, for $s\ge0$, define
\begin{equation}\label{eq:bellman-multipliers}
 \mathcal M_1(s):=
 -\beta\frac{\tfrac89+8(1-\kappa)^2}
             {\tfrac89+8(1-\kappa)^2+\alpha_{\mathrm{FS}}^{-2}s},
 \qquad
 \mathcal M_2(s):=\frac{\gamma s}{(1+\alpha_{\mathrm{FS}}^{-2}s)^{3/2}}.
\end{equation}
Evaluate the multipliers and their derivatives at $|\nabla_\theta h|^2$ unless an argument is displayed.
Adding \eqref{eq:bellman-zero-one} with $\mathcal M_1$ and
\eqref{eq:bellman-zero-two} with $\mathcal M_2$ to \eqref{eq:angular-bellman}
gives, by \eqref{eq:graph-angular-weight},
\begin{equation}\label{eq:bellman-integral-representation}
 \mathcal E_{\mathrm{ang}}[h]=\int_{\mathbb S^{d-1}}\mathcal W_h\mathcal P\,d\sigma,
\end{equation}
where
\begin{equation}\label{eq:bellman-integrand}
\begin{aligned}
 \mathcal P
 &=(1-\mathcal M_1)S+\left(\mathcal M_1-\frac1{n-1}\right)T^2+2\mathcal M_1'B_h
  +(n-2)\mathcal M_1 D_h\\
 &\quad +m_h^2\bigl[\mathcal M_2\Delta_\theta h+2\mathcal M_2'\nabla_\theta^2h(\nabla_\theta h,\nabla_\theta h)\bigr]\\
 &\quad +|\nabla_\theta h|^2\left[-\frac{n-d}{n-1}-(d-2)\mathcal M_1+(n-2)m_h^2\mathcal M_2\right].
\end{aligned}
\end{equation}

\smallskip\noindent\emph{Critical points.}
If $\nabla_\theta h=0$, then $D_h=B_h=\mathcal M_2(0)=0$,
$T=\Delta_\theta h$, and $S=\|\nabla_\theta^2h\|_{\mathrm{HS}}^2$.
By \eqref{eq:bellman-parameters} and \eqref{eq:bellman-multipliers},
\[
 1-\mathcal M_1(0)=\frac{(d-1)(n-2)}{(n-1)(d-2)},\qquad
 \mathcal M_1(0)-\frac1{n-1}=-\frac{1-\mathcal M_1(0)}{d-1}.
\]
Thus \eqref{eq:bellman-integrand} reduces to
\begin{equation}\label{eq:bellman-critical-integrand}
 \mathcal P
 =\frac{(d-1)(n-2)}{(n-1)(d-2)}
 \left\|\nabla_\theta^2h-\frac{\Delta_\theta h}{d-1}g_{\mathbb S^{d-1}}
 \right\|_{\mathrm{HS}}^2\ge0.
\end{equation}

\smallskip\noindent\emph{Points where $\nabla_\theta h\ne0$.}
Take an orthonormal frame with $e_1=\nabla_\theta h/|\nabla_\theta h|$.
Write $Y=(d-2)^{-1}\sum_{i=2}^{d-1}h_{ii}$ and set
\begin{equation}\label{eq:bellman-frame-variables}
 \xi:=\frac{\bigl(\tfrac89+8(1-\kappa)^2\bigr)m_h}
              {m_h^2+8(1-\kappa)^2-\tfrac19},\qquad 0<\frac{\xi}{m_h}<1.
\end{equation}
In this frame, \eqref{eq:graph-T}, \eqref{eq:graph-angular-hessian},
and \eqref{eq:graph-hessian-contractions} give
\begin{equation}\label{eq:bellman-frame-hessian}
\begin{gathered}
 \Delta_\theta h=h_{11}+(d-2)Y,\quad T=\frac{h_{11}}{m_h^2}+(d-2)Y, \quad S=\frac{h_{11}^2}{m_h^4}+\frac{2}{m_h^2}\sum_{i=2}^{d-1}h_{1i}^2
   +\sum_{i,j=2}^{d-1}h_{ij}^2\\
 D_h=(d-2)|\nabla_\theta h|^2Y,\qquad
 B_h=|\nabla_\theta h|^2\left((d-2)h_{11}Y-\sum_{i=2}^{d-1}h_{1i}^2\right).
\end{gathered}
\end{equation}
Define
\begin{equation}\label{eq:bellman-fjk}
 \varpi:=\left(1-\frac{\xi}{m_h}\right)(1+m_h\xi),\qquad
 k_d:=\varpi(2-\kappa\varpi)+(d-3)\left(1-\frac{\xi}{m_h}\right).
\end{equation}
By \eqref{eq:bellman-frame-variables},
\begin{equation}\label{eq:bellman-xi-product}
 m_h\xi-1=\left(8(1-\kappa)^2-\tfrac19\right)
           \left(1-\frac{\xi}{m_h}\right).
\end{equation}
The polynomial $2s+[8(1-\kappa)^2-1/9]s^2$ equals $\varpi$ at
$s=1-\xi/m_h$. Its derivative is at least $2(1-s/9)>0$ on $[0,1]$, so
\begin{equation}\label{eq:bellman-J-bound}
 0<\varpi<\tfrac{17}{9}+8(1-\kappa)^2,\qquad
 2-\kappa\varpi>2(1-\kappa)(1-2\kappa)^2+\frac{\kappa}{9}>0.
\end{equation}
In particular, $k_d>0$ by \eqref{eq:bellman-fjk}. The multipliers
in \eqref{eq:bellman-multipliers} satisfy
\begin{equation}\label{eq:bellman-multiplier-identities}
\begin{gathered}
 \mathcal M_1=-\beta\frac{\xi}{m_h},\qquad
 |\nabla_\theta h|^2\mathcal M_1'
 =\beta\frac{\xi}{m_h}\left(1-\frac{\xi}{m_h}\right),\\
 \mathcal M_1+m_h^2|\nabla_\theta h|^2\mathcal M_1'
 =-\beta(1-\varpi).
\end{gathered}
\end{equation}

Put $X:=h_{11}/m_h^2-(1-\kappa\varpi)Y$. By \eqref{eq:bellman-fjk},
\begin{equation}\label{eq:bellman-square-completion}
\begin{aligned}
 (n-2)\frac{h_{11}^2}{m_h^4}
 -2(n-2)(1-\kappa\varpi)\frac{h_{11}Y}{m_h^2}
 &+\bigl[(n-2)+(d-3)(n-d)(1-\xi/m_h)\bigr]Y^2\\
 &=(n-2)X^2+(n-d)k_dY^2.
\end{aligned}
\end{equation}
By \eqref{eq:bellman-frame-hessian}, \eqref{eq:bellman-multiplier-identities},
and \eqref{eq:bellman-square-completion}, the quadratic terms in
\eqref{eq:bellman-integrand} satisfy
\begin{equation}\label{eq:bellman-quadratic-squares}
\begin{aligned}
 (1-\mathcal M_1)S&+\left(\mathcal M_1-\frac1{n-1}\right)T^2+2\mathcal M_1'B_h\\
 &=\frac{n-2}{n-1}X^2+\frac{n-d}{n-1}k_dY^2
   +\frac{2[1+\beta(1-\varpi)]}{m_h^2}\sum_{i=2}^{d-1}h_{1i}^2\\
 &\qquad+\left(1+\frac{\beta\xi}{m_h}\right)
   \left[\sum_{i=2}^{d-1}(h_{ii}-Y)^2
         +2\sum_{2\le i<j\le d-1}h_{ij}^2\right]\\
 &\ge\frac{n-2}{n-1}X^2+\frac{n-d}{n-1}k_dY^2.
\end{aligned}
\end{equation}
Indeed, $\varpi<2/\kappa$ by \eqref{eq:bellman-J-bound}, so
\[
 1+\beta(1-\varpi)>1-\beta\,\frac{2-\kappa}{\kappa}
 =1-\frac{2-\kappa}{d-1-\kappa}\ge0.
\]

For the linear terms, \eqref{eq:bellman-multipliers} gives
$\mathcal M_2+2|\nabla_\theta h|^2\mathcal M_2'=3\gamma|\nabla_\theta h|^2/m_h^5$.
By \eqref{eq:bellman-frame-hessian}, \eqref{eq:bellman-multiplier-identities},
\begin{equation}\label{eq:bellman-linear-part}
\begin{aligned}
 (n-2)&\mathcal M_1 D_h
   +m_h^2\bigl[\mathcal M_2\Delta_\theta h+2\mathcal M_2'\nabla_\theta^2h(\nabla_\theta h,\nabla_\theta h)\bigr]\\
 &=m_h^2(\mathcal M_2+2|\nabla_\theta h|^2\mathcal M_2')h_{11}
   +(d-2)\bigl[m_h^2\mathcal M_2+(n-2)\mathcal M_1 |\nabla_\theta h|^2\bigr]Y\\
 &=\frac{\gamma |\nabla_\theta h|^2}{m_h}
   \Bigl\{3X+\bigl[(d+1)(1-\xi)-3\kappa\varpi\bigr]Y\Bigr\}.
\end{aligned}
\end{equation}

For the term independent of the Hessian, set
\begin{equation}\label{eq:bellman-constant-factors}
 \mathcal P_0:=\frac{(n-2)\gamma|\nabla_\theta h|^4}{m_h}>0,\qquad
 a_d:=\frac{(n-1)(d+1)}
           {(d-1)(n-2)^2\bigl(\tfrac89+8(1-\kappa)^2\bigr)}.
\end{equation}
By \eqref{eq:bellman-frame-variables} and
$m_h^2-1=\frac{n-1}{d-1}|\nabla_\theta h|^2$,
\begin{equation}\label{eq:bellman-constant-part}
\begin{aligned}
 |\nabla_\theta h|^2&\left[-\frac{n-d}{n-1}-(d-2)\mathcal M_1
                         +(n-2)m_h^2\mathcal M_2\right]\\
 &=\mathcal P_0-\frac{n-d}{n-1}|\nabla_\theta h|^2
                     \left(1-\frac{\xi}{m_h}\right)
        =\mathcal P_0(1-a_d\xi).
\end{aligned}
\end{equation}

Minimizing in $X$ and $Y$ in \eqref{eq:bellman-quadratic-squares}
and \eqref{eq:bellman-linear-part}, with the constant term
\eqref{eq:bellman-constant-part}, gives
\begin{equation}\label{eq:bellman-signs}
 \mathcal P\ge\mathcal P_0(1-a_d\xi-Q_{d-1}),\qquad
 Q_{d-1}:=\frac{9\kappa}{4(d+1)m_h}
 +\frac{\bigl[(d+1)(1-\xi)-3\kappa\varpi\bigr]^2}{4(d+1)m_h k_d}.
\end{equation}
On $\{\nabla_\theta h\ne0\}$, \eqref{eq:bellman-scalar-sign} in
Lemma~\ref{lem:angular-sign} with $m=m_h>1$ and \eqref{eq:bellman-signs}
give $\mathcal P>0$, since $\mathcal P_0>0$ by \eqref{eq:bellman-constant-factors}.
At critical points, $\mathcal P\ge0$ by \eqref{eq:bellman-critical-integrand}.
For nonconstant $h$ the first set is open and nonempty, so
$\mathcal E_{\mathrm{ang}}[h]>0$ by \eqref{eq:bellman-integral-representation}
and $\mathcal W_h>0$. For constant $h$, \eqref{eq:angular-bellman} gives equality.
\end{proof}

\subsection{Positivity of the scalar remainder}

\begin{lemma}\label{lem:angular-sign}
Let $n>d\ge3$ and $m>1$. Evaluate $\xi,\varpi,k_d$ and $Q_{d-1}$ in
\eqref{eq:bellman-frame-variables}, \eqref{eq:bellman-fjk}, and
\eqref{eq:bellman-signs} at $m_h=m$, with $\kappa$ from
\eqref{eq:bellman-parameters} and $a_d$ from \eqref{eq:bellman-constant-factors}.
Then
\begin{equation}\label{eq:bellman-scalar-sign}
 a_d\xi+Q_{d-1}<1.
\end{equation}
\end{lemma}

\begin{proof}
By \eqref{eq:bellman-parameters},
\begin{equation}\label{eq:bellman-sign-parameters}
 1-\kappa=\frac{d-2}{n-2}\in(0,1),\qquad
 -\frac19<8(1-\kappa)^2-\frac19<8.
\end{equation}

\smallskip\noindent\emph{Dimension reduction.}
Fix $\kappa$ and $m$. At $d=3$, \eqref{eq:bellman-fjk},
\eqref{eq:bellman-constant-factors}, and \eqref{eq:bellman-signs} give
\begin{equation}\label{eq:bellman-dimension-three}
\begin{gathered}
 a_3=\frac{2(1-\kappa)(2-\kappa)}{\tfrac89+8(1-\kappa)^2},\qquad k_3=\varpi(2-\kappa\varpi),\\
 Q_2=\frac{9\kappa}{16m}
     +\frac{[4(1-\xi)-3\kappa\varpi]^2}{16m k_3},
 \qquad \mathcal B:=\frac{(1-\xi)^2}{4(m-\xi)}.
\end{gathered}
\end{equation}
By \eqref{eq:bellman-frame-variables} and \eqref{eq:bellman-J-bound},
$0<\xi/m<1$ and $k_3>0$. Since $k_d=k_3+(d-3)(1-\xi/m)$, Cauchy--Schwarz gives
\begin{equation}\label{eq:bellman-dimension-cauchy}
 \frac{[(d+1)(1-\xi)-3\kappa\varpi]^2}{k_d}
 \le \frac{[4(1-\xi)-3\kappa\varpi]^2}{k_3}
      +(d-3)\frac{(1-\xi)^2}{1-\xi/m}.
\end{equation}
Equality holds in \eqref{eq:bellman-dimension-cauchy} at $d=3$.
Since $n-2=(d-2)/(1-\kappa)$ and $8(d-1)(d-2)-(d+1)^2=(d-3)(7d-5)\ge0$,
\eqref{eq:bellman-constant-factors} gives
\begin{equation}\label{eq:bellman-dimension-coefficient}
 \frac{(d+1)a_d}{4a_3}
 =\frac{(d+1)^2}{8(d-1)(d-2)}
  \frac{1+(1-\kappa)/(d-2)}{2-\kappa}\le1.
\end{equation}
By \eqref{eq:bellman-signs}, \eqref{eq:bellman-dimension-cauchy},
and \eqref{eq:bellman-dimension-coefficient},
\begin{equation}\label{eq:bellman-dimension-reduction}
 a_d\xi+Q_{d-1}
 \le\frac4{d+1}(a_3\xi+Q_2)+\frac{d-3}{d+1}\mathcal B.
\end{equation}

The identity
\[
9\kappa\varpi(2-\kappa\varpi)+[4(1-\xi)-3\kappa\varpi]^2=8(2-\kappa\varpi)(1-\xi)^2+2\kappa\varpi(2\xi+1)^2
\]
and \eqref{eq:bellman-dimension-three} give
\begin{equation}\label{eq:bellman-Q2-split}
 Q_2=U_0+\kappa U_1,\qquad
 U_0:=\frac{(1-\xi)^2}{2m\varpi},\qquad
 U_1:=\frac{(2\xi+1)^2}{8m(2-\kappa\varpi)}.
\end{equation}
By \eqref{eq:bellman-frame-variables} and \eqref{eq:bellman-sign-parameters},
\begin{equation}\label{eq:bellman-a3f-bound}
 a_3\xi
 =\frac{2(1-\kappa)(2-\kappa)m}{m^2-1/9+8(1-\kappa)^2}\le\frac{2(1-\kappa)(2-\kappa)m}{(8/9)m^2+8(1-\kappa)^2}
 \le\frac38(2-\kappa),
\end{equation}
where $(8/9)m^2+8(1-\kappa)^2\ge(16/3)(1-\kappa)m$.
As $m\downarrow1$, $\xi\to1$, so $a_3\le3(2-\kappa)/8$.
By \eqref{eq:bellman-xi-product}, \eqref{eq:bellman-fjk}, and \eqref{eq:bellman-J-bound},
\begin{equation}\label{eq:bellman-V-denominator}
\begin{aligned}
 2-\kappa\varpi
 &=\frac{2\xi}{m}
   +2(1-\kappa)\frac{\xi}{m}\left(1-\frac{\xi}{m}\right)+\left[2(1-\kappa)(1-2\kappa)^2+\frac{\kappa}{9}\right]
            \left(1-\frac{\xi}{m}\right)^2\\
 &\ge\frac{2\xi}{m}+\frac{\kappa}{9}\left(1-\frac{\xi}{m}\right)^2.
\end{aligned}
\end{equation}

\smallskip\noindent\emph{Case $\xi\le1$.}
Here $0\le1-\xi\le1-\xi/m<1$, and $m\xi>8/9$ by
\eqref{eq:bellman-xi-product} and \eqref{eq:bellman-sign-parameters}.
Thus $\varpi=(1-\xi/m)(1+m\xi)\ge(17/9)(1-\xi/m)$, and
\eqref{eq:bellman-dimension-three} and \eqref{eq:bellman-Q2-split} give
\begin{equation}\label{eq:bellman-small-f-basic}
 \mathcal B\le\frac{1-\xi/m}{4m}<1,\qquad
 U_0\le\frac9{34}(1-\xi)\le\frac{1-\xi}{3}.
\end{equation}
The identity $1-\xi-9\xi(1-4\xi)=(1-6\xi)^2+2\xi$ gives
$9\xi(1-4\xi)_+\le1-\xi$. Since $(m\xi)^{-1}\le9/8$,
\begin{equation}\label{eq:bellman-small-f-numerator}
\begin{aligned}
 \frac{(2\xi+1)^2-9\xi}{m}
 &=\frac{(1-\xi)(1-4\xi)}{m}
 \le\frac98\xi(1-\xi)(1-4\xi)_+\le\frac18\left(1-\frac{\xi}{m}\right)^2.
\end{aligned}
\end{equation}
Since $9\kappa<2(2+3\kappa)$ and $1/8<(2+3\kappa)/9$,
\eqref{eq:bellman-small-f-numerator} and \eqref{eq:bellman-V-denominator}
give
\begin{equation}\label{eq:bellman-small-f-V}
\begin{aligned}
 \frac{\kappa(2\xi+1)^2}{m}
 &\le\frac{9\kappa\xi}{m}+\frac{\kappa}{8}\left(1-\frac{\xi}{m}\right)^2\\
 &<(2+3\kappa)\left[\frac{2\xi}{m}
     +\frac{\kappa}{9}\left(1-\frac{\xi}{m}\right)^2\right]
 \le(2+3\kappa)(2-\kappa\varpi).
\end{aligned}
\end{equation}
By \eqref{eq:bellman-Q2-split} and \eqref{eq:bellman-small-f-V}, $\kappa U_1<(2+3\kappa)/8$.
Together \eqref{eq:bellman-small-f-basic} and \eqref{eq:bellman-a3f-bound} give
\[
 a_3\xi+U_0
 \le a_3\xi+\frac{1-\xi}{3}
 \le\max\left\{a_3,\frac13\right\}
 \le\frac38(2-\kappa).
\]
Consequently,
\begin{equation}\label{eq:bellman-small-f-conclusion}
 a_3\xi+Q_2<\frac{3(2-\kappa)+2+3\kappa}{8}=1,\qquad \mathcal B<1.
\end{equation}

\smallskip\noindent\emph{Case $\xi>1$.}
By \eqref{eq:bellman-frame-variables}, $1<m<8(1-\kappa)^2-1/9<8$, $0<\kappa<2/3$, and
\begin{equation}\label{eq:bellman-large-f-identity}
 0<\frac{\xi^2-1}{1-\xi/m}
 =\left(8(1-\kappa)^2-\tfrac19\right)\frac{\xi}{m}-1
 <8(1-\kappa)^2-\tfrac{10}{9}.
\end{equation}
For $1<s<8$,
\[
 \frac{(1+s)m}{m^2+s}\le\frac{1+s}{2\sqrt{s}}
 <\frac{9}{4\sqrt2}<\frac53.
\]
At $s=8(1-\kappa)^2-1/9$ this gives $1<\xi<5/3$. By
\eqref{eq:bellman-V-denominator} and \eqref{eq:bellman-Q2-split},
\begin{equation}\label{eq:bellman-large-f-V}
 U_1\le\frac{(2\xi+1)^2}{16\xi}
 =\frac{4\xi+4+\xi^{-1}}{16}
<\frac34.
\end{equation}
By $(1+\xi)^2\ge4\xi$, \eqref{eq:bellman-xi-product},
\eqref{eq:bellman-large-f-identity}, and \eqref{eq:bellman-dimension-three},
\begin{equation}\label{eq:bellman-large-f-B}
 \mathcal B
 =\frac{1-\xi/m}{4m(1+\xi)^2}
   \left(\frac{\xi^2-1}{1-\xi/m}\right)^2<\frac{[8(1-\kappa)^2-10/9](m\xi-1)}{16m\xi}
<1.
\end{equation}
By \eqref{eq:bellman-fjk} and \eqref{eq:bellman-Q2-split},
$U_0=2\mathcal B/(1+m\xi)$. Since
$m\xi(1+m\xi)-4(m\xi-1)=(m\xi-3/2)^2+7/4>0$,
\eqref{eq:bellman-large-f-B} yields
\begin{equation}\label{eq:bellman-large-f-U}
\begin{aligned}
 U_0
 <\frac{[8(1-\kappa)^2-10/9](m\xi-1)}{8m\xi(1+m\xi)}
 &<\frac{8(1-\kappa)^2-10/9}{32}
 <\frac{2-3\kappa}{8},
\end{aligned}
\end{equation}
since $0<\kappa<2/3$ and \eqref{eq:bellman-sign-parameters} holds.
By \eqref{eq:bellman-a3f-bound}, \eqref{eq:bellman-Q2-split},
\eqref{eq:bellman-large-f-V}, and \eqref{eq:bellman-large-f-U},
\begin{equation}\label{eq:bellman-large-f-conclusion}
 a_3\xi+Q_2<\frac{3(2-\kappa)+2-3\kappa+6\kappa}{8}=1,\qquad
 \mathcal B<1
\end{equation}
by \eqref{eq:bellman-large-f-B}. The two cases
\eqref{eq:bellman-small-f-conclusion} and \eqref{eq:bellman-large-f-conclusion}
give \eqref{eq:bellman-scalar-sign} by \eqref{eq:bellman-dimension-reduction}.
\end{proof}

\section{Proofs of the main results}
\label{sec:cone-conclusion}

\subsection{Proof of Theorem~\ref{thm:cone-isop}}

\begin{proof}[Proof of Theorem~\ref{thm:cone-isop}]
The case $\mu(E)=0$ follows from $P_\alpha(E;\Omega)\ge0$. Fix $V>0$
and choose $R$ by $\mu(B_R(0))=V$ as in \eqref{eq:cone-centered-ball}.

At $\alpha=\alpha_{\mathrm{FS}}$, Lemma~\ref{lem:geometric-compactness}
gives a minimizer of volume $V$. For any such $E$,
Lemma~\ref{lem:cone-cap-symmetrization} gives $\mu(E^{\#})=V$ and
\begin{equation}\label{eq:cone-symmetral-minimizer}
 I_{\alpha_{\mathrm{FS}}}(V)
 \le P_{\alpha_{\mathrm{FS}}}(E^{\#};\Omega)
 \le P_{\alpha_{\mathrm{FS}}}(E;\Omega)
 =I_{\alpha_{\mathrm{FS}}}(V).
\end{equation}
Thus $E^{\#}$ also minimizes perimeter. By
Proposition~\ref{prop:cap-strict-starshaped}, $(E^{\#})^{(1)}$
contains the origin in its interior and its boundary is a smooth positive
radial graph. This graph is stationary by \eqref{eq:cone-symmetral-minimizer}.
Its radius is constant by Proposition~\ref{prop:stationary-graphs} and equals
$R$ by \eqref{eq:cone-centered-ball}. Hence $(E^{\#})^{(1)}=B_R(0)$ and
\begin{equation}\label{eq:cone-critical-profile}
 I_{\alpha_{\mathrm{FS}}}(V)
 =P_{\alpha_{\mathrm{FS}}}(E^{\#};\Omega)
 =P_{\alpha_{\mathrm{FS}}}(B_R(0);\Omega).
\end{equation}

Cap symmetrization preserves slice measures by \eqref{eq:cap-spherical-slices}
and \eqref{eq:cap-opening-angle}, so
\[
 |E_r|=|(E^{\#})_r|
 =\omega_{d-1}\mathbf1_{(0,R)}(r)
 \qquad\text{for a.e. }r>0.
\]
Consequently,
\[
 \mu(E\triangle B_R(0))
 =\int_0^R r^{n-1}\bigl(\omega_{d-1}-|E_r|\bigr)\,dr
  +\int_R^\infty r^{n-1}|E_r|\,dr
 =0.
\]
Since $E$ was arbitrary,
\begin{equation}\label{eq:cone-critical-minimizers}
 \mu(E)=V,\quad P_{\alpha_{\mathrm{FS}}}(E;\Omega)=I_{\alpha_{\mathrm{FS}}}(V)
 \quad\Longrightarrow\quad \mu(E\triangle B_R(0))=0.
\end{equation}

For $0<\alpha\le\alpha_{\mathrm{FS}}$, the integrands in
\eqref{eq:polar-perimeter} satisfy the pointwise comparison
\begin{equation}\label{eq:perimeter-normal-comparison}
 \nu_r^2+\alpha^{-2}|\nu_T|^2
 \ge \nu_r^2+\alpha_{\mathrm{FS}}^{-2}|\nu_T|^2.
\end{equation}
By \eqref{eq:perimeter-normal-comparison}, \eqref{eq:cone-critical-profile},
and \eqref{eq:cone-ball-profile}, every admissible $E$ of volume $V$ satisfies
\begin{equation}\label{eq:cone-perimeter-comparison}
\begin{aligned}
 P_\alpha(E;\Omega)
 \ge P_{\alpha_{\mathrm{FS}}}(E;\Omega)
 \ge I_{\alpha_{\mathrm{FS}}}(V)
 &=P_{\alpha_{\mathrm{FS}}}(B_R(0);\Omega)\\
 &=P_\alpha(B_R(0);\Omega)=\omega_{d-1}^{1/n}(nV)^{(n-1)/n}.
\end{aligned}
\end{equation}
This proves \eqref{eq:cone-isop}. Equality in
\eqref{eq:cone-perimeter-comparison} forces $E$ to minimize at
$\alpha_{\mathrm{FS}}$, so $\mu(E\triangle B_R(0))=0$ by
\eqref{eq:cone-critical-minimizers}. Conversely, any set agreeing
a.e.\ with $B_R(0)$ has the same perimeter.
\end{proof}

\subsection{Proofs of Corollaries~\ref{cor:two-power} and~\ref{cor:stargraph-sphere}}

\begin{proof}[Proof of Corollary~\ref{cor:two-power}]
By \eqref{eq:cone-two-power-densities},
\begin{equation}\label{eq:two-power-cone-parameters}
 \alpha=L+1-K,\qquad
 n=\frac{L+d}{\alpha},\qquad
 n-1=\frac{K+d-1}{\alpha}.
\end{equation}
Since $K\ge0$, \eqref{eq:two-power-range} gives
\[
 0<\alpha\le\frac{d-1}{K+d-1}\le1,\qquad
 \alpha^2(n-1)=\alpha(K+d-1)\le d-1,
\]
and hence
\[
 0<\alpha\le\sqrt{\frac{d-1}{n-1}}=\alpha_{\mathrm{FS}},\qquad
 \alpha(n-d)=K+(d-1)(1-\alpha)\ge0.
\]
Moreover,
\[
 n=d
 \quad\Longleftrightarrow\quad K=0,\ \alpha=1
 \quad\Longleftrightarrow\quad (K,L)=(0,0).
\]
In this case \eqref{eq:two-power-isop} is Euclidean isoperimetry, with
equality at positive volume precisely for balls of arbitrary center up
to null sets \cite[Theorem~14.1]{M12}.

For $n>d$, take $E=\Phi(\widetilde E\setminus\{0\})$ as in
\eqref{eq:cone-transform}. By \eqref{eq:cone-pullback-volume-perimeter},
\[
 \mu(E)=\alpha\operatorname{Vol}_L(\widetilde E),\qquad
 P_\alpha(E;\Omega)=\operatorname{Per}_K(\widetilde E).
\]
By Theorem~\ref{thm:cone-isop} and \eqref{eq:two-power-cone-parameters},
\[
\begin{aligned}
 \operatorname{Per}_K(\widetilde E)
 =P_\alpha(E;\Omega)
 &\ge\omega_{d-1}^{1/n}\bigl(n\mu(E)\bigr)^{(n-1)/n}\\
 &=\omega_{d-1}^{(L+1-K)/(L+d)}
   \bigl((L+d)\operatorname{Vol}_L(\widetilde E)\bigr)^{(K+d-1)/(L+d)}.
\end{aligned}
\]
For every $R>0$,
\[
 \mu(E\triangle B_R(0))
 =\alpha\operatorname{Vol}_L
   \bigl(\widetilde E\triangle B_{R^{1/\alpha}}(0)\bigr).
\]
At positive weighted volume, the equality case of Theorem~\ref{thm:cone-isop}
therefore reads
\[
\begin{aligned}
 \text{equality in \eqref{eq:two-power-isop}}
 &\quad\Longleftrightarrow\quad
 \mu(E\triangle B_R(0))=0
 \quad\text{for some }R>0\\
 &\quad\Longleftrightarrow\quad
 |\widetilde E\triangle B_s(0)|=0
 \quad\text{for some }s>0.
\end{aligned}
\]
This completes the proof.
\end{proof}

\begin{remark}
The sufficient condition in \cite[Theorem~5.3]{ABCMP17},
\[
 \frac{1}{L+d}\ge \frac{1}{K+d-1}-\frac{(d-1)^2}{d(K+d-1)^3},
\]
becomes $\alpha^2\le\alpha_0^2$ in \eqref{eq:known-geometric-symmetry}
by \eqref{eq:two-power-cone-parameters}. For $K\ge L+1$ and $L+d>0$,
see \cite[Theorem~5.1]{ABCMP17}.
\end{remark}

\begin{proof}[Proof of Corollary~\ref{cor:stargraph-sphere}]
With $h=(\log f)/(n-1)$ and $\alpha_{\mathrm{FS}}^{-2}=(n-1)/(d-1)$,
\eqref{eq:def-Pv} and \eqref{eq:def-Vv} give
\[
\begin{aligned}
 \frac{P_{\alpha_{\mathrm{FS}}}(E[f^{1/(n-1)}];\Omega)}{\omega_{d-1}}
 &=\int_{\mathbb S^{d-1}}
   \sqrt{f^2+\frac{|\nabla_\theta f|^2}{(d-1)(n-1)}}\,d\sigma,\\
 \frac{n\mu(E[f^{1/(n-1)}])}{\omega_{d-1}}
 &=\int_{\mathbb S^{d-1}} f^{n/(n-1)}\,d\sigma.
\end{aligned}
\]
Thus \eqref{eq:cone-isop} gives \eqref{eq:stargraph-sphere}.
By Theorem~\ref{thm:cone-isop}, equality holds precisely when some $R>0$ satisfies
\[
 \frac{n}{\omega_{d-1}}\mu(E[f^{1/(n-1)}]\triangle B_R(0))
 =\int_{\mathbb S^{d-1}}
   \bigl|f^{n/(n-1)}-R^n\bigr|\,d\sigma=0.
\]
By continuity, this is equivalent to $f\equiv R^{n-1}$.
\end{proof}

\subsection{The cone energy space and the sharp radial quotient}
\label{sec:ckn-application}

For the cone transform $w$ of $v\in C_c^\infty(\mathbb R^d)$, the cutoff
in \eqref{eq:cone-radial-cutoff} satisfies
\[
 \int_\Omega|\nabla_\alpha((1-\chi_\varepsilon)w)|^p\,d\mu
 \le C\int_{B_{2\varepsilon}(0)}|\nabla_\alpha w|^p\,d\mu
    +C\|w\|_{L^\infty(\Omega,d\mu)}^p\varepsilon^{n-p}\to0
 \qquad(\varepsilon\downarrow0)
\]
since $n\ge d>p$. Thus \eqref{eq:cone-energy-mass} extends the cone
transform to an isomorphism from $\mathcal D_a^{1,p}$ onto
$\mathcal D^{1,p}_\alpha(\Omega,d\mu)$, with
\[
 \|v\|_{\mathcal D_a^{1,p}}
 =\alpha^{(p-1)/p}\|\nabla_\alpha w\|_{L^p(\Omega,d\mu)}.
\]
By \eqref{eq:ckn-norm},
$\mathcal D^{1,p}_\alpha(\Omega,d\mu)\hookrightarrow L^q(\Omega,d\mu)\cap W^{1,p}_{\mathrm{loc}}(\Omega)$.

Let $S_{\mathrm{rad}}$ denote the infimum of $\mathcal Q_{\alpha,p}$ in
\eqref{eq:cone-quotient-factor} over nonzero radial functions in
$\mathcal D^{1,p}_\alpha(\Omega,d\mu)$.
By \cite[Definition~2.4 and Theorem~2.2(2)]{HK12} and \cite[Theorem~1.4]{M14},
\begin{equation}\label{eq:sharp-radial-quotient}
 S_{\mathrm{rad}}
 =n\left(\frac{n-p}{p-1}\right)^{p-1}
   \left[\frac{\omega_{d-1}}{p'}\,
   \operatorname{B}\!\left(\frac n{p'},\frac np\right)\right]^{p/n}.
\end{equation}
The nonzero nonnegative radial minimizers are exactly
\begin{equation}\label{eq:radial-cone-profiles}
 W_{A,B}(r)=A\bigl(1+Br^{p'}\bigr)^{-(n-p)/p},\qquad A,B>0.
\end{equation}

\subsection{Weighted P\'olya--Szeg\H{o} inequality and equality cases}

Assume \eqref{eq:ckn-parameter-range} with $a<b<a+1$ and
$0<\alpha\le\alpha_{\mathrm{FS}}$, so $n>d\ge2$. For every measurable
$E\subset\Omega$ of finite weighted volume and perimeter,
\begin{equation}\label{eq:rearrangement-isoperimetric-hypothesis}
 P_\alpha(E;\Omega)\ge\omega_{d-1}^{1/n}\bigl(n\mu(E)\bigr)^{(n-1)/n},
\end{equation}
with equality at positive volume precisely for centered balls up to null
sets. For $d\ge3$ this is Theorem~\ref{thm:cone-isop}. For $d=2$,
\eqref{eq:cone-two-power-densities} gives $K=\alpha(n-1)-1$ and
$L=\alpha n-2$. Since $n>2$ and $\alpha\le\alpha_{\mathrm{FS}}$,
\[
 L-K+1=\alpha>0,\qquad
 2K-L=\alpha(n-2)>0,\qquad
 K^2-L(K+1)=1-\alpha^2(n-1)\ge0.
\]
These inequalities exclude $(L,K)=(0,0)$, so \cite[Theorem~1.1]{M23}
gives \eqref{eq:rearrangement-isoperimetric-hypothesis} and its equality case.

For measurable $u$ with $\mu(\{|u|>t\})<\infty$ for all $t>0$, let $u^*$
be the centered decreasing rearrangement of $|u|$ with respect to $\mu$.
Write $u^*(r)$ for its radial value and $\partial_r u^*$ for its weak
radial derivative when it exists. By equimeasurability,
\begin{equation}\label{eq:cone-rearrangement-notation}
\begin{gathered}
 R_u(t):=\left(\frac{n\mu(\{|u|>t\})}{\omega_{d-1}}\right)^{1/n},\\
 \{u^*>t\}=B_{R_u(t)}(0)\quad\text{up to null sets for every }t>0.
\end{gathered}
\end{equation}
Thus $\|u^*\|_{L^q(\Omega,d\mu)}=\|u\|_{L^q(\Omega,d\mu)}$.
By \eqref{eq:cone-metric-measure} and \eqref{eq:cone-centered-ball},
\begin{equation}\label{eq:ps-model}
 \mu(B_r(0))=\int_0^r\omega_{d-1}s^{n-1}\,ds,\qquad
 \frac{d}{dr}\mu(B_r(0))=\omega_{d-1}r^{n-1}
 =P_\alpha(B_r(0);\Omega).
\end{equation}

We use the rearrangement argument of \cite{Aub76,Tal76,BBMP99}.
Alvino et al.~\cite[Theorem~8.1 and Remark~8.2]{ABCMP17} proved weighted
P\'olya--Szeg\H{o} in their parameter range and related their isoperimetric
conjecture to equality of the full and radial CKN constants. Here
\eqref{eq:rearrangement-isoperimetric-hypothesis} covers
$\alpha^2(n-1)\le d-1$. Above this threshold, symmetry breaks by
\cite[Theorem~1.1]{CM13}. The equality criterion of Nobili and Violo
\cite[Theorem~4.2]{NV25} identifies almost every positive superlevel set
of a nonnegative optimizer with a centered ball, up to null sets.

\begin{lemma}[Weighted P\'olya--Szeg\H{o} inequality]
\label{lem:cone-rearrangement}
Assume \eqref{eq:ckn-parameter-range} with $a<b<a+1$ and
$0<\alpha\le\alpha_{\mathrm{FS}}$. For $w\in\mathcal D^{1,p}_\alpha(\Omega,d\mu)$,
the rearrangement $w^*$ belongs to the radial subspace of
$\mathcal D^{1,p}_\alpha(\Omega,d\mu)$ and satisfies
\begin{equation}\label{eq:weak-ps}
 \int_0^\infty |\partial_r w^*(r)|^p\,\omega_{d-1}r^{n-1}\,dr
 \le \int_\Omega|\nabla_\alpha w|^p\,d\mu.
\end{equation}
If $w\ge0$, $w^*(r)=W_{A,B}(r)$ a.e.\ for some $A,B>0$, and equality
holds in \eqref{eq:weak-ps}, then $w=w^*$ a.e.\ on $\Omega$.
\end{lemma}

\begin{proof}
In the completed cone $\Omega\cup\{0\}$, the perimeter relative to
an open $A\subset\Omega$ is $P_\alpha(E;A)$, and the minimal weak upper
gradient on $\Omega$ is $|\nabla_\alpha\cdot|$.

Let $0\le u\in W^{1,p}_{\mathrm{loc}}(\Omega)$ have
$\nabla_\alpha u\in L^p(\Omega,d\mu)$ and $\mu(\{u>t\})<\infty$ for $t>0$.
For $0<\varepsilon<M$, set $T_{\varepsilon,M}(s)=\min\{(s-\varepsilon)_+,M-\varepsilon\}$. Then
\[
\begin{aligned}
 &\|T_{\varepsilon,M}(u)\|_{L^p(\Omega,d\mu)}^p
 \le (M-\varepsilon)^p\mu(\{u>\varepsilon\})<\infty,\\
 &\nabla_\alpha T_{\varepsilon,M}(u)
 =\mathbf1_{\{\varepsilon<u<M\}}\nabla_\alpha u,\qquad
 \bigl(T_{\varepsilon,M}(u)\bigr)^*=T_{\varepsilon,M}(u^*).
\end{aligned}
\]
By \eqref{eq:ps-model} and \eqref{eq:rearrangement-isoperimetric-hypothesis},
\cite[Theorem~1.1]{NV25} applies to $T_{\varepsilon,M}(u)$ and gives
\begin{equation}\label{eq:ps-truncated-energy}
 \int_0^\infty \bigl|\partial_r[T_{\varepsilon,M}(u^*(r))]\bigr|^p\,\omega_{d-1}r^{n-1}\,dr
 \le \int_{\{\varepsilon<u<M\}}|\nabla_\alpha u|^p\,d\mu.
\end{equation}
On each compact interval $J\Subset(0,\infty)$, $u^*$ is bounded and
$\inf_{r\in J}\omega_{d-1}r^{n-1}>0$. Thus \eqref{eq:ps-truncated-energy} gives
\[
 \sup_{0<\varepsilon<M}
   \|T_{\varepsilon,M}(u^*)\|_{W^{1,p}(J)}<\infty,\qquad
 T_{\varepsilon,M}(u^*)\to u^*
   \quad\text{in }L^p(J)
   \qquad(\varepsilon\downarrow0,\ M\uparrow\infty).
\]
Weak compactness gives $u^*\in W^{1,p}_{\mathrm{loc}}(0,\infty)$.
By the chain rule, \eqref{eq:ps-truncated-energy} reads
\begin{equation}\label{eq:ps-band}
 \int_{\{r>0:\,\varepsilon<u^*(r)<M\}}|\partial_r u^*(r)|^p\,\omega_{d-1}r^{n-1}\,dr
 \le \int_{\{\varepsilon<u<M\}}|\nabla_\alpha u|^p\,d\mu.
\end{equation}
Let $\varepsilon\downarrow0$ and $M\uparrow\infty$ in \eqref{eq:ps-band}.
Since $\partial_r u^*=0$ a.e.\ on $\{u^*=0\}$, monotone convergence gives
\begin{equation}\label{eq:ps-general}
 \int_0^\infty |\partial_r u^*(r)|^p\,\omega_{d-1}r^{n-1}\,dr
 \le \int_\Omega|\nabla_\alpha u|^p\,d\mu.
\end{equation}
For $u=|w|$, we have $u^*=w^*$ by \eqref{eq:cone-rearrangement-notation},
$\mu(\{|w|>t\})\le t^{-q}\|w\|_{L^q(\Omega,d\mu)}^q$ for $t>0$,
and $|\nabla_\alpha u|\le|\nabla_\alpha w|$ a.e. Thus \eqref{eq:ps-general}
gives \eqref{eq:weak-ps}.

Also $w^*\in L^q(\Omega,d\mu)$ by equimeasurability. On $A_{\ell,2\ell}$
from \eqref{eq:cone-annuli}, H\"older's inequality give
\begin{equation}\label{eq:ps-cutoff}
\begin{aligned}
 \ell^{-p}\int_{A_{\ell,2\ell}}|w^*|^p\,d\mu
 &\le \ell^{-p}\mu(A_{\ell,2\ell})^{1-p/q}
       \left(\int_{A_{\ell,2\ell}}|w^*|^q\,d\mu\right)^{p/q}\\
 &\le C\left(\int_{A_{\ell,2\ell}}|w^*|^q\,d\mu\right)^{p/q}
 \to 0
 \qquad(\ell\downarrow0\ \text{or}\ \ell\uparrow\infty).
\end{aligned}
\end{equation}
By \eqref{eq:weak-ps} and \eqref{eq:ps-cutoff}, cutting off at zero and
infinity and smoothing on compact annuli approximates $w^*$ in energy
by radial $C_c^\infty(\Omega)$ functions.

Now suppose $w\ge0$, $w^*=W_{A,B}$ a.e., and equality holds in \eqref{eq:weak-ps}.
Fix $0<\varepsilon<A$. Strict decrease of $W_{A,B}$ gives
$R_w(\varepsilon)=W_{A,B}^{-1}(\varepsilon)$ by \eqref{eq:cone-rearrangement-notation}.
Since $\nabla_\alpha w=0$ a.e.\ on $\{w=\varepsilon\}$,
\begin{equation}\label{eq:ps-split}
 \int_0^\infty |W_{A,B}'|^p\,\omega_{d-1}r^{n-1}\,dr
 =\int_{\{w>\varepsilon\}}|\nabla_\alpha w|^p\,d\mu
  +\int_{\{w<\varepsilon\}}|\nabla_\alpha w|^p\,d\mu .
\end{equation}
The rearrangements of $(w-\varepsilon)_+$ and $\min\{w,\varepsilon\}$
are $(W_{A,B}-\varepsilon)_+$ and $\min\{W_{A,B},\varepsilon\}$,
respectively. By \eqref{eq:ps-general},
\begin{equation}\label{eq:ps-components}
\begin{aligned}
 \int_0^{R_w(\varepsilon)} |W_{A,B}'|^p\,\omega_{d-1}r^{n-1}\,dr
 &\le \int_{\{w>\varepsilon\}}|\nabla_\alpha w|^p\,d\mu,\\
 \int_{R_w(\varepsilon)}^\infty |W_{A,B}'|^p\,\omega_{d-1}r^{n-1}\,dr
 &\le \int_{\{w<\varepsilon\}}|\nabla_\alpha w|^p\,d\mu .
\end{aligned}
\end{equation}
The inequalities in \eqref{eq:ps-components} sum to the equality
\eqref{eq:ps-split}, so equality holds in each. In particular,
\begin{equation}\label{eq:ps-upper}
 0<\int_0^{R_w(\varepsilon)} |W_{A,B}'|^p\,\omega_{d-1}r^{n-1}\,dr
 =\int_\Omega|\nabla_\alpha (w-\varepsilon)_+|^p\,d\mu<\infty .
\end{equation}

By equimeasurability with $W_{A,B}$, $0\le w\le A$ a.e. Thus
\[
 \int_\Omega (w-\varepsilon)_+^p\,d\mu
 \le (A-\varepsilon)^p\mu(\{w>\varepsilon\})<\infty.
\]
This bound and \eqref{eq:ps-upper} give $(w-\varepsilon)_+\in W^{1,p}(\Omega,d\mu)$.
On compact annuli, the metric and density are smooth. For a.e.\ $s>0$,
the perimeter of $\{(w-\varepsilon)_+>s\}$ is concentrated on its reduced
boundary, where the approximate value equals $s$ at points of approximate
continuity. Thus \eqref{eq:cone-coarea} gives Assumption~(II) of
\cite[Section~3.2]{NV25} by exhaustion. Since
$\{(w-\varepsilon)_+>t-\varepsilon\}=\{w>t\}$ for $t>\varepsilon$,
\eqref{eq:ps-upper} and \cite[Theorem~4.2, (4.4)--(4.5)]{NV25} give
\begin{equation}\label{eq:ps-level}
 P_\alpha(\{w>t\};\Omega)
 =P_\alpha(B_{R_w(t)}(0);\Omega)
 \qquad\text{for a.e. }t\in(\varepsilon,A).
\end{equation}
Take $\varepsilon=A2^{-j}$, $j\ge1$. By \eqref{eq:ps-level} and the
equality case of \eqref{eq:rearrangement-isoperimetric-hypothesis},
\[
 \mu\bigl(\{w>t\}\triangle B_{R_w(t)}(0)\bigr)=0
 \qquad\text{for a.e. }t\in(0,A).
\]
By \eqref{eq:cone-rearrangement-notation} and Tonelli's theorem,
\[
 \int_\Omega|w-w^*|\,d\mu
 =\int_0^A\mu\bigl(\{w>t\}\triangle \{w^*>t\}\bigr)\,dt=0,
\]
so $w=w^*$ a.e.\ on $\Omega$.
\end{proof}

\subsection{Proof of Theorem~\ref{thm:main-fs}}

\begin{proof}[Proof of Theorem~\ref{thm:main-fs}]
\emph{Symmetry.}
By \eqref{eq:FS} and \eqref{eq:general-p-fs-curve},
\[
 b\ge b_{\mathrm{FS}}^{(p)}(a)
 \quad\Longleftrightarrow\quad
 \alpha^2(n-1)\le d-1
 \quad\Longleftrightarrow\quad
 0<\alpha\le\alpha_{\mathrm{FS}}.
\]
Let $w$ be the cone transform \eqref{eq:cone-transform} of $0\ne v\in\mathcal D_a^{1,p}$.
Lemma~\ref{lem:cone-rearrangement}, \eqref{eq:sharp-radial-quotient},
and equimeasurability give
\begin{equation}\label{eq:sharp-rearranged-energy}
\begin{aligned}
 \int_\Omega|\nabla_\alpha w|^p\,d\mu
 &\ge \int_0^\infty |\partial_r w^*(r)|^p\,\omega_{d-1}r^{n-1}\,dr\\
 &\ge S_{\mathrm{rad}}\|w^*\|_{L^q(\Omega,d\mu)}^p
  = S_{\mathrm{rad}}\|w\|_{L^q(\Omega,d\mu)}^p .
\end{aligned}
\end{equation}
The profiles \eqref{eq:radial-cone-profiles} attain equality in
\eqref{eq:sharp-rearranged-energy}. By \eqref{eq:cone-quotient-factor}
and $p/q=1-p/n$, we have
$C_{a,b}=C_{a,b}^{\star}=\alpha^{-(n-1)/n}S_{\mathrm{rad}}^{-1/p}$,
which gives \eqref{eq:main-constant} by \eqref{eq:sharp-radial-quotient}.

For a nonnegative optimizer $v$, equality in the second bound of
\eqref{eq:sharp-rearranged-energy} gives $w^*=W_{A,B}$ a.e.\ by
\eqref{eq:radial-cone-profiles}. Equality in the first gives $w=w^*$
a.e.\ by Lemma~\ref{lem:cone-rearrangement}. Thus
\[
 v(x)
 =W_{A,B}(|x|^\alpha)
  =A\bigl(1+B|x|^{\alpha p'}\bigr)^{-(n-p)/p}
  =A\,v_\star(\ell x)\quad\text{a.e.},\quad
 \ell=B^{(p-1)/(p\alpha)}.
\]
This proves \eqref{eq:main-optimizers}.

\smallskip\noindent
\emph{Symmetry breaking.}
Suppose $b<b_{\mathrm{FS}}^{(p)}(a)$, so $\alpha^2(n-1)>d-1$.
The gradient weight exponent in \cite[Theorem~1.1]{CM13} is $-ap$.
By \eqref{eq:parameters},
\[
 -ap>p-d,\qquad p<q=\frac{np}{n-p}<p^*,\qquad
 d-p-ap=\alpha(n-p),\qquad
 d-\frac{q}{p}(d-p-ap)=bq.
\]
The weighted Hardy inequality \cite[(1.5)]{M14} gives
\[
 (a_c-a)^p\int_{\mathbb R^d}|x|^{-(a+1)p}|v|^p\,dx
 \le \int_{\mathbb R^d}|x|^{-ap}|\nabla v|^p\,dx.
\]
Conversely, 
on $A_{s,2s}$, $|x|<2s$ gives
\[
 s^{-p}\int_{A_{s,2s}}|x|^{-ap}|v|^p\,dx
 \le 2^p\int_{A_{s,2s}}|x|^{-(a+1)p}|v|^p\,dx
 \to 0
 \qquad(s\downarrow0\ \text{or}\ s\uparrow\infty).
\]
Since $\frac1{q-p}-\frac1{q+p'}=\frac{(n-p)^2}{p^2(n-1)}$, the instability
criterion reads
\[
 \frac1{d-1}\left(\frac{d-p-ap}{p}\right)^2
 >\frac1{q-p}-\frac1{q+p'}
 \quad\Longleftrightarrow\quad
 \alpha^2(n-1)>d-1.
\]
By \cite[Theorem~1.1]{CM13}, no radial function minimizes
\eqref{eq:original-quotient}. The radial infimum is attained by
\eqref{eq:radial-cone-profiles}, hence $C_{a,b}>C_{a,b}^{\star}$.
\end{proof}

\noindent\textbf{Conflict of interest:} Authors state no conflict of interest.

\noindent\textbf{Data Availability Statement:} Data sharing is not applicable to this article as no
datasets were generated or analysed during the current study.

\noindent\textbf{Acknowledgement.}
G-D.~Li was supported by NSFC (No.~12561019).
The authors acknowledge the use of AI tools. All mathematical arguments
and proofs in the final manuscript were checked and written by the authors.

\providecommand{\href}[2]{#2}
\providecommand{\arxiv}[1]{\href{http://arxiv.org/abs/#1}{arXiv:#1}}
\providecommand{\url}[1]{\texttt{#1}}
\providecommand{\urlprefix}{DOI }
\raggedbottom

\end{document}